# Rotation excursion algorithm with learning

Sheng-Xue He

(Business School, University of Shanghai for Science and Technology, Shanghai 200093, China)

**Abstract:** We introduce a novel heuristic algorithm named the Rotation Excursion Algorithm with Learning (REAL) designed for general-purpose optimization. REAL draws inspiration from the construction mechanism inherent in CEC optimization suites, integrating three fundamental operations with a natural growth rule to address optimization tasks. The initial operation involves rotating the current feasible solutions within the search space to generate and evaluate new solutions. The excursion operation aims to relocate current feasible solutions closer to historically superior solutions stored in a list known as the "list of visible spots." The third operation involves perturbing solutions generated by the preceding operations within their respective neighborhoods. The rotation operation is geared toward comprehensive and random exploration of the entire search space, while the excursion operation exploits known information to refine current solutions. Perturbation operation functions as a form of neighborhood search to further enhance solution quality. The natural growth rule dynamically adjusts REAL's balance between exploration and exploitation throughout the entire search process. To validate the efficacy of the proposed algorithm, we apply it to address a diverse set of 67 problems, encompassing 29 benchmark optimization problems, 30 test problems from CEC 2014, one from CEC 2022, and seven engineering problems. Numerical experiments demonstrate the superior performance of REAL when compared to various other heuristics.

**Key words**: numerical optimization, heuristic algorithm, learning algorithm, swarm intelligence

## 1. Introduction

Complex optimization problems frequently arise in real-world scenarios, assuming intricate forms that pose formidable challenges for classical mathematical optimization techniques [1]. Traditional methods often rely on derivative information or specific features of functions. However, in cases where such information is unavailable or insufficient, heuristic algorithms emerge as valuable solutions. Heuristic algorithms function as versatile "black-box" tools, demanding minimal or no prior knowledge about the optimization problem. Their successful application across diverse fields over the past few decades underscores their effectiveness in tackling these intricate challenges.

Nevertheless, as highlighted by the "No Free Lunch Theorems for optimization" [2], no single general-purpose algorithm can excel in all cases. This theorem explicitly demonstrates that the performance gain of an algorithm in one problem class inevitably involves a trade-off in performance on others. Recently, Piotrowski, Napiorkowski and Piotrowska [3] have shown that the choice of different test problems leads to varying performances among algorithms. This further validates the conclusion of the No Free Lunch Theorem. Both theoretical insights and practical experiences emphasize the ongoing demand for novel heuristic algorithms in practice.

Motivated by the above observation, we propose a novel heuristic algorithm named the Rotation Excursion Algorithm with Learning (REAL) in this paper, aiming to provide practitioners with a powerful optimization tool. The foundational concepts of REAL stem from an examination of the generational process observed in the complex test functions within the CEC series [4-7]. The first of the three primary operations of REAL, the rotation operation, unfolds as follows: a feasible solution to an optimization problem undergoes rotation in a normalized search space after being transformed from the original space with a varied rotation center. The rotated solution is then re-transformed back into the original search space. The second core operation, the excursion operation, involves moving a solution towards a historical solution with outstanding performance. This operation is more likely to occur as the search progresses, guided by an increasing learning efficiency. The third operation, perturbation operation, introduces perturbations to a solution within its neighborhood, with a gradually reduced amplitude. To maintain a balance between exploitation and exploration capacities in REAL, a learning efficiency is generated based on the natural growth rule. This learning efficiency determines the likelihood of conducting an excursion operation and reduces the perturbation amplitude at each iteration of REAL. We will apply the resulting REAL to 67 benchmarks to showcase its efficiency and effectiveness.

The remainder of this paper is organized as follows: Section 2 will introduce related heuristic algorithms classified into different groups. Section 3 will elucidate the implementational procedure of REAL, first by introducing its core operations and then presenting its implementational process. Section 4 will validate the new algorithm by implementing it into various benchmark optimization problems. Section 5 will summarize the main contents and highlight potential research directions. The Appendix will present partial details of the utilized optimization problems.

**2. Literature review**

In the subsequent section, we will present a curated overview of heuristic algorithms, organized chronologically. Given space limitations, we will focus on key algorithms, excluding discussions on most variant versions. These heuristics will be systematically classified into three distinct groups. The initial group encompasses algorithms inspired by biological behaviors, excluding those derived from human social activities. The second group encompasses algorithms rooted in human social activities, while the third group includes heuristic algorithms inspired by physical phenomena and mathematical concepts.

Many algorithms in the initial group have roots dating back to earlier periods. Holland [8] introduced the Genetic Algorithm (GA), a pioneering contribution. Subsequently, Kennedy and Eberhart [9] developed the Particle Swarm Optimization (PSO) algorithm, while Storn and Price [10] introduced Differential Evolution (DE) specifically designed for global optimization in continuous spaces. The PSO algorithm was further advanced by He, Prempain and Wu [11] to address mechanical design optimization problems. Mezura-Montes and Coello [12] innovated a multimembered Evolution Strategy (ES) tailored for global, nonlinear optimization problems.

In the realm of bio-inspired algorithms, Karaboğa [13] combined the foraging behavior of honey bees to formulate the Artificial Bee Colony (ABC) algorithm. Dorigo, Birattari and Stutzle [14] extensively surveyed applications of Ant Colony Optimization (ACO). Simon [15] introduced Biogeography-Based Optimization (BBO), a novel approach inspired by biogeographical distribution principles. Drawing inspiration from the obligate brood parasitic behavior of certain cuckoos, Yang and Deb [16] devised the Cuckoo Search (CS) algorithm. Oftadeh, Mahjoob and Shariatpanahi [17] proposed the Hunting Search (HS) meta-heuristic algorithm, drawing inspiration from the cooperative hunting habits of species.

Further contributions include the development of the Firefly Algorithm (FA) and the Bat Algorithm (BA) for nonlinear design problems by researchers [18, 19]. Yang [20] conducted a comprehensive comparison of various heuristic algorithms proposed earlier. Over the course of long-term development, these early heuristic algorithms, especially GA, PSO, DE, ES, and ACO, have become widely utilized in addressing a diverse array of optimization problems.

In the past decade, there has been a proliferation of heuristic algorithms within the first group. Yang [21] introduced the Flower Pollination Algorithm (FPA), drawing inspiration from the intricate pollination process of flowers. Yang and Hossein Gandomi [22] devised the Bat Algorithm (BA) to address engineering optimization challenges. During this period, Mirjalili and collaborators contributed significantly to the field, unveiling a series of heuristic algorithms. These include the Grey Wolf Optimizer (GWO) [23], the Moth-Flame Optimization (MFO) [24] algorithm, the Whale Optimization Algorithm (WOA) [25], the Salp Swarm Algorithm (SSA) [26], and the Multiobjective Salp Swarm Algorithm (MSSA) [26]. Azizyan, Miarnaeimi, Rashki and Shabakhty [27] introduced the Flying Squirrel Optimizer (FSO). Heidari, Mirjalili, Faris, Aljarah, Mafarja and Chen [28] innovated the Harris Hawks Optimizer (HHO), a population-based optimization paradigm inspired by the cooperative hunting procedures of Harris hawks. Abdollahzadeh, Gharehchopogh and Mirjalili [29] developed the African Vultures Optimization Algorithm (AVOA), mimicking the foraging and navigational methods of African vultures. In a recent contribution, Zhao, Wang and Mirjalili [30] unveiled the Artificial Hummingbird Algorithm (AHA), simulating the flight expertise and foraging techniques of hummingbirds. Hashim and Hussien [31] proposed the Snake Optimizer, a heuristic algorithm inspired by the mating behavior of snakes, offering a novel approach to address various optimization problems.

In addition to the previously mentioned algorithms inspired by biological behaviors, a substantial body of research has focused on emulating human behaviors to create efficient algorithms. Ray and Liew [32] proposed an optimization algorithm grounded in the information exchange occurring within and between societies. Lee and Geem [33] developed the Harmony Search (HS) algorithm, incorporating the concept of musical harmonies to optimize problems with continuous variables. Rao, Savsani and Vakharia [34] introduced Teaching-Learning-Based Optimization (TLBO) specifically tailored for mechanical design optimization. He [35] introduced the Medalist Learning Algorithm (MLA), a novel swarm heuristic algorithm inspired by the learning behavior observed in groups. He [35], along with He and Cui [36], applied MLA to truss structural optimization with natural frequency constraints and the configuration optimization of trusses. Building on MLA, He and Cui [37] and He and Cui [38] modified the original algorithm to propose the Multiscale Model Learning Algorithm (MMLA), applying it to various engineering optimization problems and the search for equilibrium states in multi-tiered supply chain networks. These algorithms collectively fall within the second group mentioned earlier in this subsection. Due to their foundational ideas derived from observations of human society, these algorithms are generally straightforward to comprehend and implement in practice.

In addition to the aforementioned heuristic algorithms, many others have emerged from the observation and simulation of physical phenomena or the application of mathematical concepts. Kirkpatrick, Gelatt and Vecchi [39] pioneered the Simulated Annealing algorithm (SA), applying the concept of annealing in solid structures to optimize complex systems. Formato [40] introduced the Central Force Optimization algorithm (CFO), drawing inspiration from gravitational kinematics. Rashedi, Nezamabadi-pour and Saryazdi [41] proposed the Gravitational Search Algorithm (GSA) based on the Newtonian law of gravity and mass interactions. Kaveh and Talatahari [42] developed the Charged System Search (CSS), grounded in the Coulomb law and the Newtonian laws of mechanics. Eskandar, Sadollah, Bahreininejad and Hamdi [43] introduced the Water Cycle Algorithm (WCA), modeled on how streams and rivers navigate to the sea. Moghaddam, Moghaddam and Cheriet [44] presented the Curved Space Optimization (CSO), leveraging the concept of space-time curvature. Hatamlou [45] devised the Black Hole Algorithm (BHA) and applied it to the clustering problem. Sadollah, Bahreininejad, Eskandar and Hamdi [46] proposed the Mine Blast Algorithm (MBA) based on the concept of mine bomb explosion. Salimi [47] introduced the Stochastic Fractal Search algorithm (SFS), rooted in growth according to a mathematical fractal concept. Savsani and Savsani [48] presented the Passing Vehicle Search algorithm (PVS), founded on the mathematics of vehicles passing on a two-lane highway. Kaveh and Dadras [49] developed the Thermal Exchange Optimization algorithm (TEO). Faramarzi, Heidarinejad, Stephens and Mirjalili [50] introduced the Equilibrium Optimizer algorithm (EO), inspired by the use of balance models to estimate both dynamic and equilibrium states. Zhao, Wang and Zhang [51] created the Artificial Ecosystem-based Optimization algorithm (AEO), based on the flow of energy in an ecosystem. Most of the mentioned algorithms have demonstrated effectiveness in addressing specific optimization problems. The REAL, to be presented in this paper, belongs to this group and leverages the construction mechanism of the CEC test problems.

Recent years have witnessed active research in the field of heuristic algorithms, with numerous novel proposals making significant contributions across various domains. For instance, Ma, Wei, Tian, Cheng and Zhang [52] introduced a multi-stage evolutionary algorithm tailored for multi-objective optimization. Yildiz, Pholdee, Bureerat, Yildiz and Sait [53] improved the classical grasshopper optimization algorithm by incorporating an elite opposition-based learning regime. Zhang [54] devised an Elite Archives-driven Particle Swarm Optimization (EAPSO) to enhance the global searching ability of PSO. Ghasemi, Zare, Zahedi, Akbari, Mirjalili and Abualigah [55] introduced a new geological-inspired heuristic algorithm called Geyser-inspired Algorithm (GEA). Rezaei, Safavi, Abd Elaziz and Mirjalili [56] presented the Geometric Mean Optimizer (GMO) designed for optimizing well-known engineering problems. For a comprehensive review of existing heuristic algorithms, valuable resources include Mohammed Aarif, Sivakumar, Caffiyar, Hashim, Hashim and Rahman [57], Yang [58], and Korani and Mouhoub [59].

### 3. Rotation Excursion Algorithm with Learning (REAL)

In this section, we will commence by elucidating the foundational concepts of the REAL algorithm. Subsequently, a detailed exposition will be provided on the three primary operations and the incorporation of learning efficiency. Finally, the implementation procedure of the REAL will be presented, complete with its algorithm pseudocode and a comprehensive flow chart.

### 3.1. Fundamental ideas

The foundational concepts of REAL stem from observing the construction process of various CEC test problems [4-7]. Researchers typically manipulate the feasible search space by rotating it, shifting the global optimum, and rescaling the variables of a standard benchmark optimization problem to create complex and challenging test problems. These operations aim to disrupt the potentially regular search space and alter the position of the global optimum, rendering the resulting problems more difficult to solve.

Inspired by the effectiveness of these operations, we leverage them to design a heuristic algorithm. The rotation operation facilitates effective exploration, enabling the heuristic to traverse a broader range of the search space in a balanced manner. Shifting the position of the global optimum provides a means for the heuristic to gradually approach a potential optimum, while rescaling the variables enhances the effectiveness of the rotation operation.

In REAL, the initial step involves transforming a solution from the original search space to a standard search space, defining the feasible range of each dimension within the interval [-1,1]. The transformed solution is then rotated by multiplying it with an orthogonal rotation matrix, generating a new solution in the standard search space. Subsequently, the rotated solution is transformed back from the normalized standard search space to the original search space—a transformation akin to rescaling. In many heuristic algorithms, such as PSO [9], MLA[35] and MMLA[38], historical solutions guide the search direction. In REAL, the algorithm is inclined to move current solutions closer to historical best ones as the search progresses—an operation termed the excursion operation. The decision to conduct an excursion operation is determined by REAL based on a natural growth rule.

### 3.2. The key operations

Within this subsection, we will systematically present the rotation operation, followed by the excursion operation, the perturbation operation, and the generation and application of learning efficiency.

### 3.2.1. Rotation operation

Consider a point, also referred to as a spot, within the feasible search space—a solution to the optimization problem. To facilitate rotation, the initial step involves transforming this point into a normalized standard search space, restricting each dimension's element to the interval [-1, 1]. Subsequently, the vector linked to the transformed spot undergoes multiplication with a randomly generated orthogonal rotation matrix, yielding a new point within the normalized search space. Finally, it is imperative to revert this newly acquired point to the original search space, completing the process and generating the final point of interest.

To begin, let us explore the process in one dimension. Assuming orthogonal rotation matrices of appropriate sizes have been pre-generated. Let $M_{Rot}$ be the set of these orthogonal rotation matrices and $M \in M_{Rot}$ is one of them. For illustration, let's focus on the $i$th element of the point $x$. Represented as $x_i$, this element belongs to the feasible range $[\alpha_i, \beta_i]$ within the original search space. The rotation center is denoted as $x_C$, with $x_{C,i}$ specifically referring to the $i$th element of $x_C$.

Let $z_i \in [-1,1]$ be the corresponding element of $x_i$ in the normalized standard search space after transformation. $z_i$ can be obtained using the following two equations:

$$r_i = max\{\beta_i - x_{C,i}, x_{C,i} - \alpha_i\}, \qquad (1)$$
$$z_i = (x_i - x_{C,i})/r_i. \qquad (2)$$

Eq. (1) is to obtain the radius of $i$th dimension with $x_C$ as the center.

Let $z$ be the vector consisting of the above elements obtained by transformation. To obtain the renewed $z$ with an orthogonal rotation matrix $M$ randomly chosen from $M_{Rot}$, the rotation operation is directly applied to $z$ as follows:

$$z := Mz. \qquad (3)$$

Let $z_i$ still represent the $i$th element of the above obtained vector $z$. We need to transform it back into the original search space as follows:

$$x_i = \Pi_{[\alpha_i, \beta_i]}(z_i r_i + x_{C,i}). \qquad (4)$$

In eq. (4), the function $\Pi_{[\alpha_i, \beta_i]}(\cdot)$ is the projection operator on interval $[\alpha_i, \beta_i]$ defined as follows:

$$\Pi_{[\alpha_i, \beta_i]}(x_i) = \begin{cases} \alpha_i, & if\ x_i < \alpha_i; \\ x_i, & if\ \alpha_i \leq x_i \leq \beta_i; \\ \beta_i, & if\ x_i > \beta_i. \end{cases} \qquad (5)$$

Let $I_r$ be the diagonal matrix with all the diagonal elements equal to $r$ and all the other elements equal to zero. The above operation will be applied to all the elements of a given point $x$.

The whole operation process can be concisely expressed as follows:
$$x := \Pi_{[\alpha,\beta]}(I_r MI_r^{-1}(x - x_C) + x_C). \tag{6}$$
In eq. (6), α and β are the lower and upper boundary vectors of feasible variables in the original search space. Note that in this paper we always assume that all the vector in consideration are column vectors if not otherwise specified.

### 3.2.2. Excursion operation

In contrast to the rotation operation, the excursion operation is relatively straightforward. Let $X_{VS}$ be a list of predetermined size $n_{VS}$ storing the historically best solutions updated up to the present moment. Only solutions with superior objective function values can be added to $X_{VS}$, and when the capacity of $X_{VS}$ is reached upon adding a new solution, the one with the poorest performance in this list is replaced. Let $x_{vs}$ be a randomly chosen solution from $X_{VS}$. If the size of $X_{VS}$ is substantial, for instance, equal to 20, a roulette wheel selection mechanism may be employed here to choose a solution from $X_{VS}$ based on the associated objective function values of the solutions in $X_{VS}$. The excursion operation is executed as follows for a given point $x$:
$$x := x + r_{Ex}(x_{vs} - x), \tag{7}$$
where $r_{Ex} \in [0,1]$ is a given parameter controlling the amplitude of the excursion towards the target $x_{vs}$. Clearly, the larger the value of $r_{Ex}$, the closer the newly obtained $x$ is to $x_{vs}$. In this paper, we assign 0.5 to $r_{Ex}$ if not otherwise specified.

### 3.2.3. Perturbation operation

The perturbation operation constitutes a form of neighborhood search. Let $r_P \in [0,1]$ represent the perturbation rate, where $r_P$ is a predetermined parameter. The current amplitude of perturbation, denoted as $L_A$, undergoes modifications as the search progresses, and the method for adjusting $L_A$ will be elucidated in the subsequent subsection. Consider $Ran$ as a random number in the range [0, 1] following a uniform distribution, and $R_{norm}$ as a random number obtained from the standard normal distribution $N(0,1)$. The perturbation operation can be implemented as follows for each element of a given point $x$:
$$x_i := \begin{cases} x_i, & \text{if } Ran > r_P; \\ \Pi_{[\alpha_i,\beta_i]}(x_i + 0.3 R_{norm} L_A), & \text{otherwise.} \end{cases} \tag{8}$$
Since the standard normal distribution has a standard deviation of 1, we control the potential range of the resulting $R_{norm} L_A$ by multiplying it by 0.3. This adjustment ensures that the value of $0.3 R_{norm} L_A$ typically falls within the interval $[-3L_A, 3L_A]$. Referring to equation (8), we observe that perturbation is only carried out on the $i$th element of $x$ when $Ran \leq r_P$. For a given point $x$, REAL evaluates the perturbation operation described by equation (8) for each of its elements individually.

### 3.2.4. Generation and utilization of learning efficiency

Learning efficiency is employed to harmonize the exploration and exploitation capacities of REAL across three dimensions. First, let's examine how to generate learning efficiency associated with a learning stage based on the natural growth rule. Here, a learning stage corresponds to an iteration step of REAL. Let $T$ denote the maximum iteration number. Consider $t$ as a representative iteration of REAL, where $t \in \{1,2,\ldots,T\}$. Let $\rho(t)$ represent the learning efficiency associated with iteration $t$. The calculation of $\rho(t)$ is as follows:
$$\rho(t) = \frac{1}{1 + e^{\frac{2\gamma}{T}\left(\frac{T}{2}-t\right)}}, \tag{9}$$
where $\gamma$ is a given positive constant referred to as the radius of the truncated range of the natural growth function, specifically the Sigmoid function. According to equation (9), as $t$ increases, $\rho(t)$ initiates at a low value, gradually ascends, then enters a rapidly increasing phase, stabilizes at a high level for a duration, and ultimately converges to 1. Figure 1 illustrates the curves of $\rho(t)$ for $T=100$ with varying values of $\gamma$ set at 4, 6, and 10. Figure 3 demonstrates that for a given $T$, a larger value of $\gamma$ results in a steeper rising part of the curve. The value of $\gamma$ governs the shape of the learning efficiency curve, where a larger $\gamma$ extends the initial and final phases associated with low and high efficiencies, respectively, but narrows the middle phase of accelerating learning efficiency.

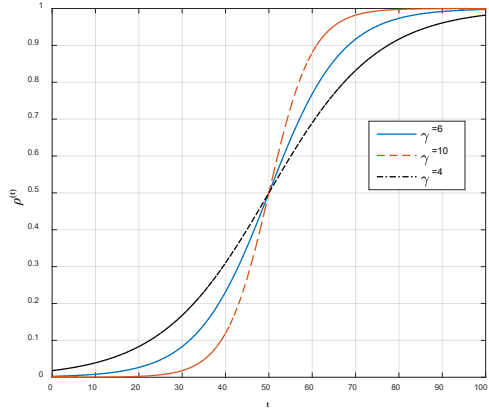

Fig. 1. Learning efficiency over the learning period $T = 100$.

As mentioned earlier, learning efficiency is employed to equilibrate the exploitation and exploration capacities of REAL. Initially, we utilize learning efficiency to modify the value of $L_A(t)$. Let $L_A^0$ and $L_A^T$ represent the initial and final amplitudes of perturbation. As the implementation of REAL progresses, $L_A(t)$ should gradually decrease from $L_A^0$ to $L_A^T$ to shift the focus of neighborhood search to a smaller area in later stages. To achieve this objective, the adjustment of $L_A(t)$ is as follows:

$$L_A(t) = L_A^0 - (L_A^0 - L_A^T)\rho(t). \qquad (10)$$

Secondly, learning efficiency can determine whether to execute an excursion operation during the search process. After implementing the rotation operation on a point $x$, we assess whether to perform an excursion operation for it. Initially, a random number $Ran$ in the range $[0,1]$, following a uniform distribution, is generated. If $Ran$ is less than $\rho(t)$, the excursion operation is applied to $x$; otherwise, the excursion operation is skipped for $x$ in the current stage $t$.

Thirdly, learning efficiency is employed to decide whether to change the rotation center $x_C$ from the center $\bar{x}_C$ of the original search space to a historical best solution $x_{vs}$ in $X_{vs}$. The center $\bar{x}_C$ of the original search space is given by $\bar{x}_C = \frac{1}{2}(\alpha + \beta)$. The approach is akin to judging whether to conduct an excursion operation. A random number $Ran$ is generated for comparison with $\rho(t)$. If $Ran \geq \rho(t)$, $\bar{x}_C$ becomes the rotation center; otherwise, a historical best solution $x_{vs}$ randomly chosen from $X_{vs}$ assumes the role of the rotation center.

### 3.3. The implementation procedure

The REAL's implementation procedure, detailed in the preceding sections, can be succinctly captured through its pseudocode, showcased in Table 1, and visually represented by the corresponding flow chart in Figure 2. Initially, a set of alternative solutions to the problem is randomly generated from the original feasible search space. The core functionality of the REAL involves applying the rotation, possible excursion, and perturbation operations to all alternative solutions, repeating this process at most $T$ times. Termination occurs either when the maximum iteration number is reached or when the best performance among the solutions in the historical best solutions list equals the worst performance of those historical solutions. Ultimately, the best solution is output as the result. Table 1 denotes $x_{vs}^{best}$ and $x_{vs}^{worst}$ as the solutions with the best and worst performances, respectively, in the list $X_{VS}$. The performance of a solution $x$ is represented by $F(x)$, with the convention that smaller values of $F(x)$ correspond to better solution performance.

If the REAL does not conclude its process prior to reaching the maximum number of iterations, it will approximately assess the objective function value about $2.5T$ times for a typical $x$ within $X$. This estimation includes the excursion operation, typically performed $0.5T$ times after $T$ iterations, following the symmetric probability curve of learning efficiency throughout the entire search period. The maximal Number of objective Function Estimations (NFE) for the REAL remains strictly below $3T$ times for each $x$ in $X$. Assuming the size of $X$ is denoted as $n_X$, the maximal NFE of the REAL is calculated as $3Tn_X$, with the NFE typically hovering around $2.5Tn_X$.

The time-intensive aspect of the REAL primarily resides in its rotation operation, attributed to matrix multiplication. Assuming the variable $x$ possesses $n$ independent elements, the rotation operation in one iteration requires approximately $n^2 n_X$ basic algebraic operations. This insight leads to the realization that the time complexity of the REAL is $O(Tn^2 n_X)$.

It is imperative to note that generating rotation matrices in advance is crucial for optimizing computational efficiency in the REAL. These matrices will be precomputed following the method outlined by Salomon [60]. Preparing rotation matrices in advance proves advantageous in saving substantial computational time during the implementation of the REAL. Delay can be markedly reduced by generating matrices proactively, contrasting with the inefficiency associated with on-the-fly matrix generation during the REAL's execution.

**Table 1**

The pseudocode of the REAL.

---

**Initialize parameters:**
  Assign values to $\bar{x}_C$, $L_A^0$, $L_A^T$, $\gamma$, $r_P$, $r_{Ex}$, $M_{Rot}$, and $T$.
**Generate the set $X$ of initial solutions:**
  Assume that the size of $X$ is given. For any $x \in X$, its $i$th element $x_i$ is chosen from $[\alpha_i, \beta_i]$ randomly according to a uniform distribution.
**for** $t = 1:T$
  Compute $\rho(t)$ and $L_A(t)$ with eq. (9) and eq. (10), respectively.
  **for** each $x \in X$, do the following
    Choose a matrix $M \in M_{Rot}$ randomly.
    Apply rotation operation to $x$ according to eq. (6).
    **if** $Ran < \rho(t)$
      Choose a point $x_{vs}$ from $X_{VS}$ randomly.
      Apply the excursion operation to $x$ according to eq. (7).
    **end**
    Apply perturbation operation to $x$ according to eq. (8).
    Try adding the newly obtained $x$ into $X_{VS}$.
  **end**
  **if** $F(x_{vs}^{best}) = F(x_{vs}^{worst})$, **break**;
**end**
**Output:** $F(x_{vs}^{best})$ and $x_{vs}^{best}$.

---

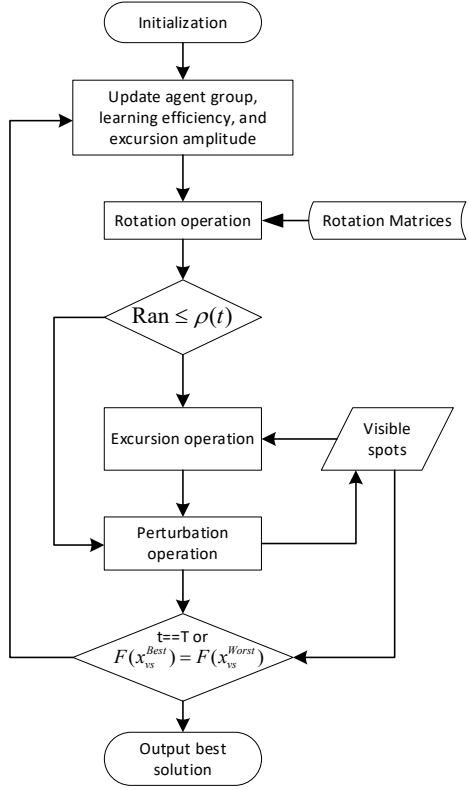

Fig. 2. The flow chart of the REAL.

**4. Numerical experiments**

In this section, we will evaluate the effectiveness of REAL by implementing it to solve 67 benchmark optimization problems. REAL is coded in JAVA and executed on Apache NetBeans IDE 12.6. The computer utilized features an Intel(R) Core (TM) i7-1065G7 CPU with 16.0 GB of RAM.

In this paper, a straightforward approach is employed to address discrete variables. When the rotation and excursion operations are executed, the discrete variable is treated as a continuous variable. Following these operations, the resulting variable is substituted with its nearest feasible discrete value. For instance, consider 4.1 and 7.7 as two consecutive discrete values within a feasible range of a one-dimensional discrete variable. If a continuous value $x_i$ is obtained such that $x_i \in$ [4.1,7.7], we compare $x_i - 4.1$ and $7.7 - x_i$. If $x_i - 4.1$ is less than $7.7 - x_i$, then $x_i := 4.1$; otherwise, $x_i := 7.7$. When considering the perturbation operation, the discrete value $x_i$ will assume either of its two neighbors with equal probability if the condition $Ran \leq r_P$ is met.

Let $n_{Rot}$ represent the number of matrices in the set $M_{Rot}$. $n_{VS}$ is the maximum number of historical best solutions stored in the set $X_{VS}$. As defined earlier, $n_X$ is the number of solutions in the group $X$. The benchmark values assigned to parameters are as follows: $n_X$=30, γ=6, T=1000, $n_{Rot}$=20, $n_{VS}$=10, $r_{A0} = 0.5$, $L_A^T = 1E - 40$, $r_P$=0.1, and $r_{Ex} = 0.5$. To ensure a fair comparison with other algorithms, we will adjust the values of $n_X$ and $T$ to make the expected number of function evaluations (NFE) (approximately $2.5n_X T$) about the same as the others. The values of the other parameters will remain unchanged unless otherwise specified.

In subsection 4.1, we will employ 29 well-known benchmark optimization problems to evaluate the effectiveness of REAL. Subsection 4.2 involves solving 30 problems from the CEC 2014 test suite using REAL and other algorithms. Subsection 4.3 will showcase the application of REAL to seven renowned engineering optimization problems, demonstrating its effectiveness in real-world engineering scenarios. Lastly, in subsection 4.4, we will use a problem from the CEC 2022 test suite to investigate the impact of varied parameters on the performance of REAL.

Please be aware that we will execute the program 30 times to address each problem in every scenario. The statistical results are derived from the average values obtained from these 30 independent runs.

**4.1. Standard optimization problems**

In this subsection, we will employ REAL to solve 29 standard benchmark optimization problems, categorized into seven unimodal problems, six high-dimensional multimodal problems, ten multimodal problems with fixed dimensions, and six composition problems. These benchmark functions have been extensively utilized in previous research to evaluate various heuristic algorithms. Notably, algorithms such as MMLA[37], EO [50], GA [8], PSO [9], GSA [41], GWO [23], SSA [26], CMA-ES [61], SHADE [62], and LSHADE [63] have been previously applied to these problems. Detailed descriptions of these benchmark problems can be found in Appendix A.

Table 2 summarizes the primary results obtained using REAL, providing a comparative analysis with results from nine other algorithms, namely MMLA, EO, PSO, GWO, GSA, SSA, CMA-ES, SHADE, and LSHADE.

Among the 29 problems, REAL outperforms other algorithms in 18 problems based on the final best average objective values (Ave). For comparison, we utilize the Friedman mean rank and assign ranks to all algorithms, as shown in Table 2. The top three algorithms are REAL, MMLA, and EO, respectively. Notably, REAL excels in solving unimodal problems and performs admirably in multimodal and composition problems. MMLA, EO, and LSHADE outperform others in 14, 7, and 6 problems, respectively. The data also reveal that each algorithm considered here demonstrates the best performance when handling at least one problem. Additionally, alongside average objective values, Table 2 presents standard deviations (Std), consistent with the performance in Ave. In other words, REAL demonstrates relatively consistent performance. While it does not outperform all other algorithms across all problems, its promising performance aligns with the "No Free Lunch Theorem"[2] and underscores its theoretical and practical value.

**Table 2**

Optimization results and comparison for functions.

| | Function | | REAL | MMLA | EO | PSO | GWO | GSA | SSA | CMA-ES | SHADE | LSHADE |
|---|---|---|---|---|---|---|---|---|---|---|---|---|
| Unimodal | F1 | Ave | 5.07E−59 | 4.95E−23 | 3.32E−40 | 9.59E−06 | 6.59E−28 | 2.53E−16 | 1.58E−07 | 1.42E−18 | 1.42E−09 | 0.2237 |
| | | Std | 5.59E−59 | 4.92E−24 | 6.78E−40 | 3.35E−05 | 1.58E−28 | 9.67E−17 | 1.71E−07 | 3.13E−18 | 3.09E−09 | 0.1480 |
| | F2 | Ave | 1.63E−19 | 3.12E−10 | 7.12E−23 | 0.02560 | 7.18E−17 | 0.05565 | 2.66293 | 2.98E−07 | 0.0087 | 21.1133 |
| | | Std | 2.53E−19 | 1.97E−11 | 6.36E−23 | 0.04595 | 7.28E−17 | 0.19404 | 1.66802 | 1.7889 | 0.0213 | 9.5781 |
| | F3 | Ave | 4.73E−23 | 2.14E−10 | 8.06E−09 | 82.2687 | 3.29E−06 | 896.534 | 1709.94 | 1.59E−05 | 15.4352 | 88.7746 |
| | | Std | 8.97E−23 | 2.79E−10 | 1.60E−08 | 97.2105 | 1.61E−05 | 318.955 | 11242.3 | 2.21E−05 | 9.9489 | 47.2300 |
| | F4 | Ave | 3.17E−16 | 2.87E−10 | 5.39E−10 | 4.26128 | 5.61E−07 | 7.35487 | 11.6741 | 2.01E−06 | 0.9796 | 2.1170 |
| | | Std | 2.06E−16 | 3.22E−11 | 1.38E−09 | 0.67730 | 1.04E−06 | 1.74145 | 4.1792 | 1.25E−06 | 0.7995 | 0.4928 |
| | F5 | Ave | 20.90517 | 24.06265 | 25.32331 | 92.4310 | 26.81258 | 67.5430 | 296.125 | 36.7946 | 24.4743 | 28.8255 |
| | | Std | 0.282880 | 0.952351 | 0.169578 | 74.4794 | 0.793246 | 62.2253 | 508.863 | 33.4614 | 11.2080 | 0.8242 |
| | F6 | Ave | 0 | 0 | 8.29E−06 | 8.89E−06 | 0.816579 | 2.5E−16 | 1.80E−07 | 6.83E−19 | 5.31E−10 | 0.2489 |
| | | Std | 0 | 0 | 5.02E−06 | 9.91E−06 | 0.482126 | 1.74E−16 | 3.00E−07 | 6.71E−19 | 6.35E−10 | 0.1131 |
| | F7 | Ave | 2.25E-5 | 0.001634 | 0.001171 | 0.02724 | 0.002213 | 0.08944 | 0.1757 | 0.0275 | 0.0235 | 0.0047 |
| | | Std | 1.18E-5 | 8.85E−04 | 6.54E−04 | 0.00804 | 0.001996 | 0.04339 | 0.0629 | 0.0079 | 0.0088 | 0.0019 |

|  | | | | | | | | | | | |
|---|---|---|---|---|---|---|---|---|---|---|---|
| Multimodal (High dimensional) | F8 | Ave | -7344.81 | -11863.47 | −9016.34 | −6075.85 | −6123.1 | −2821.1 | −7455.8 | −7007.1 | −11713.1 | −3154.4 |
| | | Std | 828.3404 | 133.4013 | 595.1113 | 754.632 | 909.865 | 493.037 | 772.811 | 773.94 | 230.49 | 317.921 |
| | F9 | Ave | 38.8170 | 0 | 0 | 52.8322 | 0.31052 | 25.9684 | 58.3708 | 25.338 | 8.5332 | 67.542 |
| | | Std | 18.2076 | 0 | 0 | 16.7068 | 0.35214 | 7.47006 | 20.016 | 8.5539 | 2.1959 | 10.016 |
| | F10 | Ave | 6.48E-15 | 1.28E-12 | 8.34E−14 | 0.00501 | 1.06E−13 | 0.06208 | 2.6796 | 15.587 | 0.3957 | 0.0393 |
| | | Std | 1.62E-15 | 4.90E-13 | 2.53E−14 | 0.01257 | 2.24E−13 | 0.23628 | 0.8275 | 7.9273 | 0.5868 | 0.0151 |
| | F11 | Ave | 3.14E-16 | 4.08E-25 | 0 | 0.02381 | 0.00448 | 27.7015 | 0.0160 | 5.76E−15 | 0.0048 | 0.8948 |
| | | Std | 5.62E-16 | 6.03E-26 | 0 | 0.02870 | 0.00665 | 5.04034 | 0.0112 | 6.18E−15 | 0.0077 | 0.1078 |
| | F12 | Ave | 7.39E-14 | 4.01E-15 | 7.97E−07 | 0.02764 | 0.05343 | 1.79961 | 6.9915 | 2.87E−16 | 0.0346 | 8.18E−04 |
| | | Std | 3.25E-14 | 6.39E-16 | 7.69E−07 | 0.05399 | 0.02073 | 0.95114 | 4.4175 | 5.64E−16 | 0.0875 | 0.0010 |
| | F13 | Ave | 1.18E-12 | 2.44E-24 | 0.029295 | 0.00732 | 0.65446 | 8.89908 | 15.8757 | 3.66E−04 | 7.32E−04 | 0.0102 |
| | | Std | 6.02E-13 | 4.21E-25 | 0.035271 | 0.01050 | 0.00447 | 7.12624 | 16.1462 | 0.0020 | 0.0028 | 0.0103 |
| Multimodal (Fixed-dimensional) | F14 | Ave | 0.998004 | 0.998004 | 0.998004 | 3.84902 | 4.042493 | 5.859838 | 1.1965 | 10.237 | 0.998004 | 1.9416 |
| | | Std | 0 | 0 | 1.54E−16 | 3.24864 | 4.252799 | 3.831299 | 0.5467 | 7.5445 | 5.83E−17 | 2.9633 |
| | F15 | Ave | 3.0749E-4 | 0.000307 | 0.002398 | 0.002434 | 0.00337 | 0.003673 | 0.000886 | 0.0057 | 0.002374 | 3.00E−04 |
| | | Std | 0 | 1.37E-19 | 0.006097 | 0.006081 | 0.001647 | 0.000257 | 0.000121 | 0.0121 | 0.0061 | 1.93E−19 |
| | F16 | Ave | -1.03163 | -1.03163 | −1.03162 | −1.03162 | −1.03163 | −1.03163 | −1.03163 | −1.03162 | −1.03162 | −1.03162 |
| | | Std | 1.40E-16 | 0 | 6.04E−16 | 6.51E−16 | 2.13E−08 | 4.88E−16 | 6.13E−14 | 6.77E−16 | 6.51E−16 | 1.00E−15 |
| | F17 | Ave | 0.397887 | 0.397887 | 0.397887 | 0.397887 | 0.397889 | 0.397887 | 0.397887 | 0.397887 | 0.397887 | 0.397887 |
| | | Std | 0 | 0 | 0 | 0 | 2.13E−04 | 0 | 3.41E−14 | 0 | 3.24E−16 | 0 |
| | F18 | Ave | 3 | 3 | 3 | 3 | 3.000028 | 3 | 3 | 8.4000 | 3 | 3 |
| | | Std | 9.42E-16 | 0 | 1.56E−15 | 1.97E−15 | 4.24E−04 | 4.17E−15 | 2.20E−13 | 20.550 | 1.87E−15 | 1.25E−15 |
| | F19 | Ave | 3.86278 | -3.86278 | −3.86278 | −3.86278 | −3.86263 | −3.86278 | −3.86278 | −3.86278 | −3.86278 | −3.86278 |
| | | Std | 8.88E-16 | 0 | 2.59E−15 | 2.65E−15 | 0.00273 | 2.29E−15 | 1.47E−10 | 2.7E−15 | 2.69E−15 | 2.7E−15 |
| | F20 | Ave | -3.3219 | -3.3220 | −3.2687 | −3.26651 | −3.28654 | −3.31778 | −3.2304 | −3.2903 | −3.27047 | −3.28234 |
| | | Std | 2.63E-15 | 0 | 0.05701 | 0.06032 | 0.10556 | 0.023081 | 0.0616 | 0.0535 | 0.0599 | 0.0570 |
| | F21 | Ave | -10.1532 | -10.1532 | −8.55481 | −5.9092 | −8.7214 | −5.95512 | −9.6334 | −5.6642 | −9.2343 | −9.4735 |
| | | Std | 1.12E-15 | 0 | 2.76377 | 3.59559 | 2.6914 | 3.73707 | 1.8104 | 3.3543 | 2.4153 | 1.7626 |
| | F22 | Ave | -10.4029 | -10.4029 | −9.3353 | −7.3360 | −9.2415 | −10.4015 | −9.0295 | −8.4434 | −10.1479 | −10.2258 |
| | | Std | 2.38E-15 | 0 | 2.43834 | 3.47381 | 1.61254 | 2.01408 | 2.3911 | 3.3388 | 1.3969 | 0.9704 |
| | F23 | Ave | -10.5364 | -10.5364 | −9.63655 | −8.7482 | −10.5343 | −10.5364 | −9.0333 | −8.0750 | −10.2809 | −10.5364 |
| | | Std | 1.68E-15 | 0 | 2.38811 | 2.55743 | 0.00125 | 2.6E−15 | 2.9645 | 3.5964 | 1.3995 | 1.77E−15 |
| Composition | F24 (CF1) | Ave | 8.71E-21 | 2.41E-27 | 66.666 | 151.18 | 90.229 | 20.000 | 43.333 | 209.48 | 63.333 | 3.3333 |
| | | Std | 1.70E-22 | 7.56E-28 | 95.893 | 123.49 | 105.51 | 48.423 | 67.891 | 215.06 | 80.872 | 18.254 |
| | F25 (CF2) | Ave | 17.35 | 26.34 | 89.837 | 204.92 | 163.56 | 186.77 | 31.133 | 189.83 | 40.508 | 0.0000 |
| | | Std | 30.077 | 22.365 | 56.366 | 118.89 | 89.476 | 62.726 | 52.149 | 170.79 | 61.462 | 0.0000 |
| | F26 (CF3) | Ave | 4.61E-21 | 2.07E-28 | 161.73 | 273.73 | 210.61 | 218.55 | 235.11 | 274.20 | 139.48 | 104.29 |
| | | Std | 2.11E-22 | 1.23E-28 | 33.227 | 110.87 | 95.214 | 117.02 | 80.839 | 213.89 | 33.366 | 14.266 |
| | F27 (CF4) | Ave | 1.47 | 6.31 | 356.44 | 487.45 | 418.63 | 492.33 | 232.44 | 372.99 | 316.62 | 278.63 |
| | | Std | 5.45E-1 | 6.67E-08 | 115.66 | 151.15 | 156.16 | 99.549 | 43.643 | 152.12 | 96.752 | 7.0670 |
| | F28 (CF5) | Ave | 11.541 | 0.600 | 52.309 | 214.56 | 143.81 | 232.32 | 27.538 | 224.85 | 39.515 | 2.02E−17 |
| | | Std | 5.632 | 0.236 | 95.565 | 180.03 | 149.12 | 75.405 | 41.598 | 286.23 | 51.233 | 7.69E−17 |
| | F29 (CF6) | Ave | 119.63 | 160.00 | 768.48 | 794.50 | 837.47 | 845.47 | 628.69 | 845.26 | 684.51 | 540.23 |
| | | Std | 132.04 | 219.09 | 192.94 | 175.94 | 136.45 | 80.524 | 184.48 | 139.52 | 201.22 | 122.75 |
| Friedman mean rank | | | 2.4138 | 2.431 | 4.6724 | 7.6034 | 6.2759 | 7.069 | 6.7414 | 6.9655 | 5.2759 | 5.5517 |
| Rank | | | 1 | 2 | 3 | 10 | 6 | 9 | 7 | 8 | 4 | 5 |

## 4.2. CEC 2014

In this section, we will assess the effectiveness of REAL using the well-known CEC 2014 test suite [4]. The results obtained with REAL will be compared to those generated by seven other algorithms: AHA[30], TLBO[34], GSA, Artificial Bee Colony (ABC) [13], Covariant Matrix Adaptation with Evolution Strategy (CMA-ES) [64], Success-History Based Adaptive DE (SHADE) [65], and Salp Swarm Algorithm (SSA) [26]. The CEC 2014 test suite comprises 30 problems, including 3 unimodal functions (CF1–CF3), 13 multimodal functions (CF4–CF16), 6 hybrid functions (CF17–CF22), and 8 composition functions (CF23–CF30). To conserve space, we will focus solely on the 30-dimensional versions of these problems. All these problems share the same search range, defined by $[-100,100]^D$, where $D=30$ represents the number of dimensions. Each problem CF$i$ is associated with a theoretical optimum such that $F_i^* = 100 \times i$ for all $i \in \{1,2,\ldots,30\}$.

With the exception of the REAL results, all other data is sourced from the literature [30]. To ensure a fair comparison, we maintain a maximum number of function evaluations (NFE) that does not exceed 25,000, as per the study [30]. For REAL, we set $T = 500$ and $n_X = 20$, keeping the other parameters consistent with the previous subsection. Specific parameters used by other algorithms can be referenced in Zhao (2022) [30]. Each of the eight algorithms is run independently 30 times, and the results presented here are based on the average performance across these 30 runs.

The outcomes related to the CEC 2014 test suite are summarized in Table 3. Within this table, the best-performing results for each problem are underlined. The data in Table 3 highlights that REAL outperforms other algorithms in 12 out of the 30 problems. AHA attains the best results for 10 problems, while SHADE leads in 6 other problems. It's worth noting that, except for ABC and CMA-ES, other algorithms also excel in specific problems.

To provide a comprehensive assessment of the overall performance of these eight algorithms, we conducted the Friedman test. The final ranking of each algorithm is presented in Table 3. Notably, REAL, AHA, and SHADE secure the top three ranks. The mean rank of REAL, which is 2.2167, signifies its superior overall performance in tackling these 30 problems.

**Table 3**

Results on functions in CEC 2014.

| Fun. | Index | REAL | AHA | TLBO | GSA | ABC | CMA-ES | SHADE | SSA |
|---|---|---|---|---|---|---|---|---|---|
| CF1 | Ave | 8.92E+6 | 2.27E+7 | <u>6.81E+6</u> | 4.20E+8 | 2.22E+8 | 1.89E+8 | 9.91E+6 | 2.58E+7 |
|  | Std | 3.42E+6 | 1.22E+7 | 3.43E+6 | 8.23E+7 | 4.81E+7 | 8.33E+6 | 2.82E+6 | 1.33E+7 |
| CF2 | Ave | <u>7.94E+3</u> | 9.42E+6 | 9.85E+5 | 2.33E+10 | 3.07E+8 | 3.89E+9 | 1.01E+4 | 1.03E+4 |
|  | Std | 2.49E+3 | 6.98E+6 | 1.16E+6 | 3.26E+9 | 7.11E+7 | 2.69E+9 | 4.40E+3 | 8.60E+3 |
| CF3 | Ave | 3.13E+3 | 6.47E+3 | 4.27E+4 | 8.73E+4 | 8.00E+4 | 2.66E+4 | <u>5.08E+2</u> | 6.78E+4 |
|  | Std | 9.33E+2 | 3.59E+3 | 6.45E+3 | 6.51E+3 | 1.36E+4 | 1.85E+4 | 1.26E+2 | 1.65E+4 |
| CF4 | Ave | <u>4.94E+2</u> | 5.67E+2 | 5.63E+2 | 2.83E+3 | 7.95E+2 | 1.94E+3 | 5.01E+2 | 5.60E+2 |
|  | Std | 14.4425 | 4.77E+1 | 3.69E+1 | 3.37E+2 | 4.52E+1 | 3.37E+2 | 2.92E+1 | 4.12E+1 |
| CF5 | Ave | 519.9999 | 520.1335 | 521.0691 | <u>519.9992</u> | 520.0542 | 520.9916 | 521.0091 | 519.9999 |
|  | Std | 6.75E-6 | 9.33E-2 | 4.66E-2 | 4.90E-4 | 5.95E-2 | 8.92E-2 | 4.54E-2 | 1.11E-4 |
| CF6 | Ave | <u>6.05E+2</u> | 6.17E+2 | 6.15E+2 | 6.32E+2 | 6.32E+2 | 6.28E+2 | 6.19E+2 | 6.20E+2 |
|  | Std | 1.09E+0 | 2.68E+0 | 2.46E+0 | 1.37E+0 | 2.74E+0 | 4.93E+0 | 4.90E+0 | 3.33E+0 |
| CF7 | Ave | 700.0144 | 701.0583 | 700.7527 | 956.9402 | 703.4847 | 790.1209 | <u>700.0042</u> | 700.0133 |
|  | Std | 8.27E-3 | 7.40E-2 | 2.61E-1 | 3.35E+1 | 4.94E-1 | 3.47E+1 | 8.88E-3 | 9.59E-3 |
| CF8 | Ave | 9.23E+2 | <u>8.48E+2</u> | 8.67E+2 | 9.45E+2 | 1.02E+3 | 9.17E+2 | 9.10E+2 | 9.42E+2 |
|  | Std | 1.82E+1 | 1.75E+1 | 1.39E+1 | 1.01E+1 | 1.15E+1 | 3.69E+1 | 1.01E+1 | 3.04E+1 |
| CF9 | Ave | 1.03E+3 | 1.05E+3 | <u>1.01E+3</u> | 1.07E+3 | 1.15E+3 | 1.10E+3 | 1.06E+3 | 1.05E+3 |
|  | Std | 2.06E+1 | 3.09E+1 | 3.14E+1 | 1.44E+1 | 1.09E+1 | 4.51E+1 | 8.42E+0 | 4.04E+1 |
| CF10 | Ave | 3.85E+3 | <u>1.67E+3</u> | 5.46E+3 | 4.91E+3 | 7.48E+3 | 3.56E+3 | 5.04E+3 | 4.47E+3 |
|  | Std | 2.04E+2 | 3.96E+2 | 1.32E+3 | 3.60E+2 | 2.55E+2 | 9.45E+2 | 5.22E+2 | 6.58E+2 |
| CF11 | Ave | <u>3.75E+3</u> | 3.89E+3 | 8.42E+3 | 5.61E+3 | 8.63E+3 | 7.94E+3 | 7.43E+3 | 4.88E+3 |
|  | Std | 4.96E+2 | 5.14E+2 | 2.94E+2 | 5.07E+2 | 2.81E+2 | 8.24E+2 | 3.71E+2 | 6.86E+2 |
| CF12 | Ave | <u>1200.0223</u> | 1200.2553 | 1203.1427 | 1200.0275 | 1203.0937 | 1202.7605 | 1202.4971 | 1200.6903 |
|  | Std | 4.63E-2 | 1.45E-1 | 3.75E-1 | 1.31E-2 | 4.56E-1 | 5.85E-1 | 3.14E-1 | 3.74E-1 |
| CF13 | Ave | <u>1300.3722</u> | 1300.4547 | 1300.4339 | 1304.4934 | 1300.6252 | 1302.8761 | 1300.5450 | 1300.4914 |
|  | Std | 4.41E-2 | 1.20E-1 | 8.87E-2 | 3.90E-1 | 8.36E-2 | 4.65E-1 | 4.32E-2 | 1.08E-1 |
| CF14 | Ave | <u>1400.1978</u> | 1400.3336 | 1400.2649 | 1504.4364 | 1400.2499 | 1445.3904 | 1400.2850 | 1400.3555 |
|  | Std | 4.49E-2 | 1.51E-1 | 4.84E-2 | 1.39E+1 | 6.97E-2 | 1.48E+1 | 3.86E-2 | 1.53E-1 |
| CF15 | Ave | 1533.8271 | 1517.8706 | 1527.2681 | 6393.9443 | 1600.5196 | 1602.8043 | 1514.3350 | <u>1513.4816</u> |
|  | Std | 5.06E+0 | 9.43E+0 | 6.09E+0 | 2.61E+3 | 3.99E+1 | 2.98E+2 | 9.37E-1 | 3.89E+0 |
| CF16 | Ave | 1611.6175 | <u>1610.9303</u> | 1612.7586 | 1613.7006 | 1613.0812 | 1612.9198 | 1612.8083 | 1612.4193 |
|  | Std | 5.23E-1 | 5.95E-1 | 3.01E-1 | 3.01E-1 | 2.12E-1 | 2.93E-1 | 2.66E-1 | 4.07E-1 |
| CF17 | Ave | <u>2.37E+5</u> | 1.32E+6 | 4.29E+5 | 3.77E+7 | 5.36E+6 | 6.70E+6 | 2.88E+5 | 1.57E+6 |
|  | Std | 2.26E+5 | 9.18E+5 | 2.67E+5 | 8.06E+6 | 1.42E+6 | 4.99E+6 | 1.72E+5 | 1.41E+6 |
| CF18 | Ave | <u>1.98E+3</u> | 3.32E+3 | 2.79E+3 | 1.98E+5 | 1.30E+4 | 1.51E+6 | 1.00E+4 | 6.67E+6 |
|  | Std | 2.76E+1 | 1.70E+3 | 1.18E+3 | 1.07E+6 | 1.54E+4 | 5.49E+6 | 6.24E+3 | 4.70E+3 |
| CF19 | Ave | 1.92E+3 | 1.93E+3 | 1.92E+3 | 2.10E+3 | 1.92E+3 | 1.97E+3 | <u>1.91E+3</u> | 1.92E+3 |
|  | Std | 3.15E+0 | 2.83E+1 | 2.41E+1 | 2.36E+1 | 1.31E+0 | 3.03E+1 | 7.83E-1 | 1.59E+1 |
| CF20 | Ave | 8.07E+3 | 1.85E+4 | 2.00E+4 | 2.80E+4 | 4.25E+4 | 3.61E+4 | <u>3.47E+3</u> | 2.86E+4 |
|  | Std | 3.36E+3 | 7.19E+3 | 5.94E+3 | 1.28E+5 | 1.51E+4 | 2.41E+4 | 2.02E+3 | 1.39E+4 |
| CF21 | Ave | 9.31E+4 | 2.91E+5 | 1.71E+5 | 1.63E+7 | 1.06E+6 | 1.06E+6 | <u>4.95E+4</u> | 3.05E+5 |
|  | Std | 4.58E+4 | 2.17E+5 | 1.09E+5 | 4.06E+6 | 5.06E+5 | 1.36E+6 | 4.13E+4 | 2.27E+5 |
| CF22 | Ave | <u>2.57E+3</u> | 2.85E+3 | 2.60E+3 | 3.70E+3 | 2.89E+3 | 3.04E+3 | 2.65E+3 | 2.67E+3 |
|  | Std | 9.99E+1 | 1.94E+2 | 1.40E+2 | 5.65E+2 | 1.41E+2 | 2.64E+2 | 8.48E+1 | 2.20E+2 |
| CF23 | Ave | 2621.1883 | <u>2500.0000</u> | 2615.3530 | 2573.1800 | 2621.9003 | 2663.2842 | 2615.2475 | 2633.2248 |
|  | Std | 2.55E+0 | 0.00E+0 | 1.48E-1 | 1.33E+0 | 1.15E+2 | 2.65E-3 | 8.08E+0 |  |
| CF24 | Ave | 2600.0134 | <u>2600.0000</u> | 2600.0800 | 2620.0005 | 2648.2497 | 2636.3027 | 2624.9547 | 2639.0200 |
|  | Std | 4.64E-3 | 0.00E+0 | 1.32E-2 | 5.73E+0 | 3.37E+0 | 6.10E+0 | 9.92E-1 | 8.32E+0 |
| CF25 | Ave | <u>2700.0000</u> | <u>2700.0000</u> | 2700.0159 | 2705.6276 | 2732.6953 | 2723.8740 | 2711.2547 | 2715.0444 |
|  | Std | 0.00E+0 | 0.00E+0 | 8.70E-2 | 2.69E+0 | 3.65E+0 | 6.90E+0 | 1.29E+0 | 3.05E+0 |
| CF26 | Ave | <u>2700.4957</u> | 2763.5233 | 2713.7514 | 2793.1042 | 2707.7046 | 2703.2826 | 2700.8536 | 2700.5293 |
|  | Std | 1.34E-1 | 4.88E+1 | 3.44E+1 | 2.16E+1 | 1.96E+1 | 1.08E+0 | 6.37E-2 | 1.46E-1 |
| CF27 | Ave | 3013.2273 | <u>2900.0000</u> | 3271.3167 | 4557.5228 | 3502.8503 | 3470.9875 | 3108.8516 | 3194.4532 |
|  | Std | 1.47E+2 | 0.00E+0 | 1.25E+2 | 3.96E+2 | 1.22E+2 | 2.29E+2 | 8.07E+1 | 1.64E+2 |
| CF28 | Ave | 4164.6349 | 3000.0000 | 3938.6520 | 5215.9261 | 4357.7299 | 7445.9885 | 3929.0812 | 4653.9440 |
| | Std | 4.24E+2 | 0.00E+0 | 1.42E+2 | 8.98E+2 | 1.18E+2 | 5.62E+2 | 7.64E+1 | 5.52E+2 |
| CF29 | Ave | 9.34E+3 | 3.10E+3 | 4.69E+3 | 6.25E+6 | 3.78E+5 | 2.85E+7 | 1.26E+4 | 1.35E+4 |
| | Std | 3.57E+3 | 0.00E+0 | 9.96E+2 | 2.02E+7 | 1.97E+5 | 3.16E+7 | 2.41E+3 | 5.47E+3 |
| CF30 | Ave | 1.75E+4 | 4.45E+3 | 6.98E+3 | 2.57E+6 | 4.70E+4 | 1.01E+5 | 8.44E+3 | 6.39E+4 |
| | Std | 5.21E+3 | 3.53E+3 | 1.52E+3 | 6.84E+5 | 1.32E+4 | 6.05E+4 | 9.29E+2 | 6.13E+4 |
| Friedman mean rank | | 2.2167 | 3.0667 | 3.7167 | 6.5167 | 6.2833 | 6.35 | 3.3333 | 4.5167 |
| Rank | | 1 | 2 | 4 | 8 | 6 | 7 | 3 | 5 |

## 4.3. Engineering optimization problems

In this subsection, we will apply REAL to solve seven well-known engineering optimization problems. These problems, derived from real-life scenarios, present a significant challenge to heuristic algorithms due to their complex search spaces and multiple local minima. Table 4 provides an overview of the key characteristics of these engineering optimization problems. The "Mixed/Con/Discrete" column in Table 4 indicates the nature of the variables involved. "Discrete" signifies that all variables are discrete, "continuous" implies that all variables fall within corresponding continuous intervals, and "Mixed" indicates a combination of discrete and continuous variables. Detailed mathematical formulations of these seven engineering problems are provided in Appendix B. It's worth noting that we modified the original model of the Robot gripper optimization problem by adding three constraints to ensure result feasibility.

To facilitate meaningful comparisons with other algorithms in terms of the total Number of objective Function Evaluations (NFE), we will adjust the total number of iterations of REAL to align the resulting NFEs with those of other algorithms.

In this paper, we employ penalty terms to address constraint violations. When a constraint is breached, the absolute difference between the boundary and the actual function value resulting from the constraint is multiplied by a large positive constant and added to the objective function.

Table 5 presents a comprehensive summary of the primary results obtained for these engineering problems using REAL.

**Table 4**

The main characters of the seven engineering optimization problems.

| No. | Name of problem | Number of Variables | Number of Constraints | Mixed/Con/Discrete |
|---|---|---|---|---|
| 1 | Multiple disc clutch brake design | 5 | 8 | Discrete |
| 2 | Robot gripper | 7 | 10 | Continuous |
| 3 | Rolling element bearing | 10 | 7 | Mixed |
| 4 | Hydrodynamic thrust bearing | 4 | 7 | Continuous |
| 5 | Belleville spring | 4 | 7 | Continuous |
| 6 | Step-cone pulley | 5 | 8 | Continuous |
| 7 | Speed reducer design | 7 | 11 | Continuous |

**Table 5**

The main results from REAL for the seven engineering problems.

| No. of Problem | 1 | 2 | 3 | 4 | 5 | 6 | 7 |
|---|---|---|---|---|---|---|---|
| Best | 0.313657 | 4.97437160 | 81859.7407 | 1625.443 | 1.9806 | 18.4484 | 2994.4711 |
| Mean | 0.313657 | 5.52179623 | 81859.7175 | 1626.127 | 2.0129 | 20.5688 | 2994.4711 |
| Worst | 0.313657 | 6.19262749 | 81859.6690 | 1627.014 | 2.1625 | 23.3296 | 2994.4711 |
| Std | 0.0 | 0.37655366 | 0.0177 | 1.135 | 0.057 | 1.596 | 3.334E-7 |
| NFE | 1200 | 30000 | 16000 | 48000 | 24000 | 72000 | 32000 |

To showcase the effectiveness of REAL, we have assembled Tables 6 to 12, comparing the results obtained using REAL with those from other algorithms. For a fair comparison, the results from other algorithms were obtained from relevant literature and represent their best-performing outcomes.

Addressing the problems sequentially as listed in Table 4, various algorithms have been applied.

The first problem has been tackled by AHA [30], TLBO [34], ABC [34], EOBL-GOA [53], AVOA [29], AEO [51], HHO [28], FSO [27], PVS [48], and WCA [43]. The second problem has been approached using multiple existing heuristic methods, including an approach combining the grasshopper optimization algorithm and the Nelder–Mead algorithm (HGOANM) [66], TLBO [34], PVS [48] and a way using Force and Displacement Transmission Ratio(2014)[67]. The third problem has been explored previously employing EAPSO [54], TLBO [34], AHA [30], AVOA [29], AEO [51], HHO [28], MBA [46] and WCA [43]. The fourth problem has been optimized using AHA [30], TLBO [34], EOBL-GOA [53], and PVS [48]. The fifth problem has been addressed with the involvement of AHA [30], PVS [48], MBA [46], TLBO [34], a multi-objective optimization approach (Coello) [68], the combined genetic search technique (Gene AS I and Gene AS II) [69], and the APPROX (Griffith and Stewart's successive linear approximation) [70]. The sixth problem has been studied using PVS [48], ABC [34], TLBO [34], and an improved constrained Differential Evolution variant (rank-iMDDE) [71]. The seventh problem has been optimized with AHA [30], Geometric Mean Optimizer(GMO)[56], EAPSO [54], AEO [51], PVS [48], MBA [46], CS [72], BA [22], WCA [43], and Artificial Bee Colony (ABC) algorithm[73].

The optimal objective values for each problem are highlighted in red, along with the best means of the objective values, representing the average values of 30 independent runs, and their associated standard deviations. The comparison illustrates that REAL consistently achieves the best objective values across all seven problems, demonstrating its potential as a heuristic algorithm.

It is crucial to note that the feasibility of all results presented in Tables 6 to 12 has been thoroughly verified. Some results may exhibit minor violations of related constraints, indicated by underlined data entries in the associated tables. These violations could be attributed to two factors: firstly, the use of low-accuracy criteria in certain studies, and secondly, truncation of results in specific papers to meet formatting constraints.

**Table 6**

Results for the multiple disc clutch brake design.

|  |  | REAL | MMLA | AHA | TLBO | EOBL-GOA | AEO | HHO | FSO |
|---|---|---|---|---|---|---|---|---|---|
| $x$ | $x_1$ | 70 | 70 | 70 | 70 | 70 | 70 | 70 | 70 |
|  | $x_2$ | 90 | 90 | 90 | 90 | 90 | 90 | 90 | 90 |
|  | $x_3$ | 1 | 1 | 1 | 1 | 1 | 1 | 1 | 1 |
|  | $x_4$ | 820 | 830 | 840 | 810 | 984 | 810 | 1000 | 870 |
|  | $x_5$ | 3 | 3 | 3 | 3 | 3 | 3 | 2.313 | 3 |
| $g(x)$ | $g_1(x)$ | 0 | 0 | 0 | 0 | 0 | 0 | 0 | 0 |
|  | $g_2(x)$ | 24 | 24 | 24 | 24 | 24 | 24 | 25.5 | 24 |
|  | $g_3(x)$ | 0.919 | 0.917 | 0.916 | 0.919 | 0.902 | 0.919 | 0.901 | 0.913 |
|  | $g_4(x)$ | 9.828 | 9.826 | 9.824 | 9.830 | 9.794 | 9.830 | 9.791 | 9.818 |
|  | $g_5(x)$ | 7.895 | 7.895 | 7.895 | 7.895 | 7.895 | 7.895 | 7.895 | 7.895 |
|  | $g_6(x)$ | 0.443 | 0.619 | 1.198 | 0.702 | 2.869 | 0.702 | <u>-2.905</u> | 0.876 |
|  | $g_7(x)$ | 38.913 | 40.119 | 41.325 | 37.706 | 58.695 | 37.706 | 20.417 | 44.944 |
|  | $g_8(x)$ | 14.557 | 14.381 | 13.802 | 14.298 | 12.131 | 14.298 | 17.905 | 14.124 |
| Best |  | 0.313657 | 0.313657 | 0.313657 | 0.313657 | 0.313657 | 0.313657 | 0.259769 | 0.313657 |
| Mean |  | 0.313657 | 0.315225 | 0.321684 | 0.327166 | 0.313657 | 0.321684 | NA | 0.313942 |
| Worst |  | 0.313657 | 0.321498 | 0.333260 | 0.392071 | 0.321612 | 0.333260 | NA | 0.333705 |
| std |  | 0.0 | 0.003137 | 0.008 | 0.67 | 0.264 | 0.008154 | NA | 4.25e-6 |
| NFE |  | 1200 | 600 | 600 | 600 | 300 | 500 | NA | 400 |

**Table 7**

Results for the robot gripper design.

|  |  | REAL | HGOANM | TLBO | PVS |
|---|---|---|---|---|---|
| $x$ | $x_1$ | 150.0 | 150.000 | 150 | 150 |
|  | $x_2$ | 131.323 | 149.883 | 150 | 150 |
|  | $x_3$ | 186.550 | 200 | 200 | 200 |

|   |        | 17.481     | 0         | 0         | 0         |
|---|--------|------------|-----------|-----------|-----------|
|   | $x_4$  | 17.481     | 0         | 0         | 0         |
|   | $x_5$  | 103.224    | 150       | 150       | 150       |
|   | $x_6$  | 145.000    | 100.943   | 100       | 100       |
|   | $x_7$  | 2.383      | 2.297     | 2.340     | 2.312     |
| $g(x)$ | $g_1(x)$  | 41.564 | 49.989 | -537.514 | -545.098 |
|   | $g_2(x)$  | 8.436      | 0.011     | 587.514   | 595.098   |
|   | $g_3(x)$  | 45.558     | 49.992    | 33.650    | 43.777    |
|   | $g_4(x)$  | 4.442      | 0.008     | 16.350    | 6.223     |
|   | $g_5(x)$  | 57812.122  | 79740.197 | 80000.0   | 80000.0   |
|   | $g_6(x)$  | 2340.496   | 36.009    | 0.0       | 0.0       |
|   | $g_7(x)$  | 45.000     | 0.943     | 0.0       | 0.0       |
|   | $g_8(x)$  | 29.599     | 0.826     | 0.0       | 0.0       |
|   | $g_9(x)$  | 127.373    | 100.826   | 100.0     | 100.0     |
|   | $g_{10}(x)$ | 135.273  | 198.940   | 200.0     | 200.0     |
| Best |     | 4.97437160 | 71.05999  | 74.9999999 | 74.99999 |
| Mean |     | 5.52179623 | NA        | NA        | NA        |
| Worst |    | 6.19262749 | NA        | NA        | NA        |
| std |      | 0.37655366 | NA        | NA        | NA        |
| NFE |      | 30000      | NA        | 25000     | 25000     |

**Table 8**

Results for the rolling element bearing design.

|   |        | REAL       | MMLA       | AHA        | EAPSO      | AVOA       | AEO        | HHO        | MBA        |
|---|--------|------------|------------|------------|------------|------------|------------|------------|------------|
| $x$ | $x_1$ | 125.719 | 125.719 | 125.718 | 125.719 | 125.723 | 125.719 | 125.0 | 125.715 |
|   | $x_2$  | 21.426     | 21.426     | 21.425     | 21.426     | 21.423     | 21.426     | 21.0       | 21.423     |
|   | $x_3$  | 11.023     | 10.920     | 10.528     | 11         | 11.001     | 11.395     | 11.092     | 11.0       |
|   | $x_4$  | 0.515      | 0.515      | 0.515      | 0.515      | 0.515      | 0.515      | 0515       | 0.515      |
|   | $x_5$  | 0.515      | 0.515      | 0.515      | 0.515      | 0.515      | 0.515      | 0.515      | 0.515      |
|   | $x_6$  | 0.486      | 0.479      | 0.470      | 0.5        | 0.404      | 0.410      | 0.4        | 0.489      |
|   | $x_7$  | 0.680      | 0.667      | 0.641      | 0.644      | 0.619      | 0.638      | 0.6        | 0.628      |
|   | $x_8$  | 0.300      | 0.300      | 0.300      | 0.3        | 0.3        | 0.300      | 0.3        | 0.300      |
|   | $x_9$  | 0.099      | 0.081      | 0.095      | 0.047      | 0.069      | 0.047      | 0.050      | 0.097      |
|   | $x_{10}$ | 0.708    | 0.629      | 0.682      | 0.605      | 0.602      | 0.670      | 0.6        | 0.646      |
| $g(x)$ | $g_1(x)$ | 2.905E-9 | 3.066e-10 | 3.900E-5 | -8.657E-10 | 0.001 | 3.148E-8 | 0.106 | 5.638E-4 |
|   | $g_2(x)$  | 8.806  | 9.328      | 9.936      | 7.851      | 14.537     | 14.137     | 14.0       | 8.630      |
|   | $g_3(x)$  | 4.750  | 3.806      | 2.007      | 2.199      | 0.461      | 1.838      | 0.0        | 1.101      |
|   | $g_4(x)$  | -0.174 | -2.551     | -0.958     | -3.265     | -3.349     | -3.141     | -3.0       | -2.041     |
|   | $g_5(x)$  | 0.719  | 0.719      | 0.718      | 0.719      | 0.723      | 0.719      | 0.0        | 0.715      |
|   | $g_6(x)$  | 24.219 | 19.499     | 23.062     | 11.086     | 16.560     | 11.023     | 12.619     | 23.611     |
|   | $g_7(x)$  | 3.728  | 3.281E-7   | 2.574E-4   | 8.000E-9   | 6.300E-6   | 2.548      | 0.7        | 5.179E-4   |
|   | $g_8(x)$  | 0.0    | 3.667E-12  | 0.0        | 0.0        | 0.0        | 0.0        | 0.0        | 0          |
|   | $g_9(x)$  | 0.0    | 1.775E-11  | 1.550E-4   | 0.0        | 0.0        | 9.000E-7   | 0.0        | 0          |
| Best |     | 81859.7407 | 81859.7392 | 81812.0128 | 81859.7416 | 81843.4096 | 81859.2991 | 78897.8107 | 81843.686  |
| Mean |     | 81859.7175 | 81777.3303 | NA         | NA         | NA         | NA         | NA         | NA         |
| Worst |    | 81859.6690 | 81344.6726 | NA         | NA         | NA         | NA         | NA         | NA         |
| std |      | 0.0177     | 120.4525   | NA         | NA         | NA         | NA         | NA         | NA         |
| NFE |      | 16000      | 15000      | 15000      | 15000      | NA         | 10000      | NA         | 15100      |

**Table 9**

Results for the hydrodynamic thrush bearing design.

|   |        | REAL       | MMLA       | TLBO       | AHA        | EOBL-GOA   | PVS        |
|---|--------|------------|------------|------------|------------|------------|------------|
| $x$ | $x_1$ | 5.955  | 5.955      | 5.956      | 5.956      | 5.956      | 5.956      |
|   | $x_2$  | 5.389      | 5.389      | 5.389      | 5.389      | 5.389      | 5.389      |
|   | $x_3$  | 5.359E-6   | 5.359E-6   | 5.4E-6     | 5.359E-6   | 5.359E-6   | 1.0E-5     |
|   | $x_4$  | 2.270      | 2.269      | 2.270      | 2.270      | 2.26966338 | 2.270      |
| $g(x)$ | $g_1(x)$ | 2.891E-4 | 1.713E-8 | 1.375E-4 | 0.041 | NA | NA |
|   | $g_2(x)$  | 1.650E-4   | 4.038E-9   | 1.0E-6     | 1.811E-4   | NA         | NA         |
|   | $g_3(x)$  | 4.801E-6   | 8.96132e-7 | 2.0E-8     | 7.629E-4   | NA         | NA         |
|   | $g_4(x)$  | 3.244E-4   | 3.24363e-4 | 3.244E-4   | 3.244E-4   | NA         | NA         |
|   | $g_5(x)$  | 0.567      | 0.567      | 0.567      | 0.567      | NA         | NA         |
|   | $g_6(x)$  | 8.334E-4   | 8.334E-4   | 8.334E-4   | 9.964E-4   | NA         | NA         |
|   | $g_7(x)$  | 4.403E-6   | 2.523E-9   | 9.08E-6    | 4.646e-3   | NA         | NA         |
| Best |     | 1625.443   | 1625.443   | 1625.443   | 1625.450   | 1625.443   | 1625.444   |
| Mean |     | 1626.127   | 1644.001   | 1797.708   | 1680.781   | 1645.422   | 1832.492   |
| Worst |    | 1627.014   | 1935.610   | 2096.801   | 1850.381   | 1675.493   | NA         |
| std |      | 1.135      | 61.108     | NA         | 57         | 12.615     | NA         |
| NFE |      | 48000      | 50000      | 50000      | 50000      | 20000      | 25000      |

**Table 10**

Results for the Belleville spring design.

|   |   | REAL | MMLA | AHA | TLBO | PVS | MBA | Coello | Gene ASII |
|---|---|------|------|-----|------|-----|-----|--------|-----------|
| $x$ | $x_1$ | 0.204 | 0.204 | 0.204 | 0.204 | 0.204 | 0.204 | 0.208 | 0.210 |
|   | $x_2$ | 0.200 | 0.200 | 0.200 | 0.2 | 0.2 | 0.2 | 0.2 | 0.204 |
|   | $x_3$ | 10.014 | 10.012 | 10.030 | 10.030 | 10.030 | 10.030 | 8.751 | 9.268 |
|   | $x_4$ | 11.997 | 11.995 | 12.010 | 12.01 | 12.01 | 12.01 | 11.067 | 11.499 |
| $g(x)$ | $g_1(x)$ | 0.008 | 1.262E-5 | 2.235 | 0.506 | 2.498 | 0.243 | 2145.411 | 2127.262 |
|   | $g_2(x)$ | 0.041 | 6.285E-6 | -0.267 | -0.036 | -0.274 | -0.028 | 39.7502 | 194.223 |
|   | $g_3(x)$ | 8.449E-7 | 4.715E-10 | 0.0 | 0.0 | 0.0 | 0.0 | 0.0 | 0.004 |
|   | $g_4(x)$ | 1.596 | 1.596 | 1.596 | 1.596 | 1.596 | 1.596 | 1.592 | 1.586 |
|   | $g_5(x)$ | 0.013 | 0.015 | 5.000E-6 | 0.0 | 0.0 | 0.0 | 0.943 | 0.511 |
|   | $g_6(x)$ | 1.983 | 1.983 | 1.980 | 1.980 | 1.980 | 1.980 | 2.316 | 2.231 |
|   | $g_7(x)$ | 0.199 | 0.199 | 0.199 | 0.199 | 0.199 | 0.199 | 0.214 | 0.209 |
| Best |   | 1.9806 | 1.9807 | 1.9797 | 1.9797 | 1.9797 | 1.9797 | 2.1220 | 2.1626 |
| Mean |   | 2.0129 | 2.0437 | 1.9860209 | 1.9886513 | 1.983524 | 1.984698 | NA | NA |
| worst |   | 2.1625 | 2.1910 | 2.1041923 | 2.1149819 | NA | 2.005431 | NA | NA |
| std |   | 0.057 | 0.054 | 0.023 | 0.027 | NA | 7.78E−03 | NA | NA |
| NFE |   | 24000 | 25000 | 24000 | 24000 | 15000 | 10600 | 24000 | 24000 |

**Table 11**

Results for the step-cone pulley design.

|   |   | REAL | MMLA (Strict) | TLBO | PVS | ABC | Rank-iMDDE |
|---|---|------|---------------|------|-----|-----|------------|
| $x$ | $x_1$ | 34.585 | 36.066 | 40 | 40 | NA | 100 |
|   | $x_2$ | 47.586 | 49.625 | 54.764 | 54.764 | NA | 34.582 |
|   | $x_3$ | 63.444 | 66.162 | 73.013 | 73.013 | NA | 47.582 |
|   | $x_4$ | 76.075 | 79.332 | 88.428 | 88.428 | NA | 63.438 |
|   | $x_5$ | 99.990 | 95.983 | 85.986 | 85.986 | NA | 76.067 |
| $g(x)$ | $h_1(x)$ | 0.0 | 0.0 | 0.001 | 0.001 | NA | 0.371 |
|   | $h_2(x)$ | 0.0 | 0.0 | 9.999E-4 | 9.999E-4 | NA | 0.367 |
|   | $h_3(x)$ | 0.0 | 0.0 | 0.001 | 0.001 | NA | 0.352 |
|   | $g_1(x)$ | 0.989 | 0.988 | 0.987 | 0.986864 | NA | 0.963 |
|   | $g_2(x)$ | 0.999 | 0.998 | 0.997 | 0.997360 | NA | 0.999 |
|   | $g_3(x)$ | 1.009 | 1.009 | 1.011 | 1.010154 | NA | 1.008 |
|   | $g_4(x)$ | 1.019 | 1.019 | 1.021 | 1.020592 | NA | 1.0156 |
|   | $g_5(x)$ | 705.658 | 707.061 | 698.577 | 698.577228 | NA | 2211.360 |
|   | $g_6(x)$ | 486.668 | 487.841 | 475.827 | 475.827090 | NA | 19.148 |
|   | $g_7(x)$ | 216.927 | 217.805 | 209.037 | 209.036940 | NA | -116.522 |
|   | $g_8(x)$ | 5.804E-6 | 0.635 | -1.528 | -1.528 | NA | -204.632 |
| Best |   | 18.4484 | 19.2588 | 21.1797 | 21.1797 | 16.6347 | 36.8448 |
| Mean |   | 20.5688 | 24.7556 | 24.0114 | NA | 36.0995 | NA |
| worst |   | 23.3296 | 29.1133 | 74.0230 | NA | 145.4705 | NA |
| std |   | 1.596 | 2.187 | 0.34 | NA | 0.06 | NA |
| NFE |   | 72000 | 75000 | NA | 15000 | 15000 | NA |

**Table 12**

Results for the speed reducer design.

|   |   | REAL | MMLA | GMO | AHA | EAPSO | AEO | MBA | CS | BA | WCA |
|---|---|------|------|-----|-----|-------|-----|-----|-----|-----|-----|
| $x$ | $x_1$ | 3.5 | 3.5 | 3.5 | 3.5 | 3.5 | 3.5 | 3.5 | 3.502 | 3.5 | 3.5 |
|   | $x_2$ | 0.7 | 0.7 | 0.7 | 0.7 | 0.7 | 0.7 | 0.7 | 0.7 | 0.7 | 0.7 |
|   | $x_3$ | 17.0 | 17.0 | 17.0 | 17.0 | 17.0 | 17 | 17 | 17 | 17 | 17 |
|   | $x_4$ | 7.3 | 7.3 | 7.3 | 7.300 | 7.3 | 7.3 | 7.300 | 7.605 | 7.300 | 7.3 |
|   | $x_5$ | 7.715 | 7.715 | 7.728 | 7.715 | 7.715 | 7.7153 | 7.716 | 7.818 | 7.715 | 7.715 |
|   | $x_6$ | 3.350 | 3.350 | 3.350 | 3.350 | 3.350 | 3.3502 | 3.350 | 3.352 | 3.350 | 3.350 |
|   | $x_7$ | 5.287 | 5.287 | 5.287 | 5.287 | 5.287 | 5.287 | 5.287 | 5.288 | 5.288 | 5.287 |
| $g(x)$ | $g_1(x)$ | -0.074 | -0.074 | -0.074 | -0.074 | -0.074 | -0.074 | -0.074 | -0.074 | -0.074 | -0.739 |
|   | $g_2(x)$ | -0.198 | -0.198 | -0.198 | -0.198 | -0.198 | -0.198 | -0.198 | -0.198 | -0.198 | -0.198 |
|   | $g_3(x)$ | -0.499 | -0.499 | -0.499 | -0.499 | -0.499 | -0.499 | -0.499 | -0.435 | -0.499 | -0.499 |
|   | $g_4(x)$ | -0.905 | -0.905 | -0.904 | -0.905 | -0.905 | -0.905 | -0.905 | -0.901 | -0.905 | -0.905 |
|   | $g_5(x)$ | -1.965E-10 | -1.506E-10 | 1.313E-5 | -2.409E-8 | -3.496E-9 | 1.432E-13 | -2.930E-6 | -0.001 | 4.195E-6 | 5.965E-7 |
|   | $g_6(x)$ | -1.696E-11 | -7.861E-11 | -2.344E-5 | -1.036E-9 | 2.824E-9 | 5.918E-13 | 3.505E-7 | -4.0E-4 | -4.797E-4 | 2.637E-7 |
|   | $g_7(x)$ | -0.703 | -0.702 | -0.703 | -0.703 | -0.703 | -0.703 | -0.703 | -0.703 | -0.703 | -0.703 |
|   | $g_8(x)$ | -1.577E-14 | -7.372E-11 | 0.0 | -2.431E-8 | 0.0 | 6.8E-13 | 0.0 | -4.0E-4 | 0.0 | 0.0 |

| | | | | | | | | | | |
|---|---|---|---|---|---|---|---|---|---|---|
| | $g_9(x)$ | -0.583 | -0.583 | -0.583 | -0.796 | -0.583 | -0.796 | -0.583 | -0.583 | -0.583 | -0.5833 |
| | $g_{10}(x)$ | -0.051 | -0.051 | -0.051 | -0.051 | -0.051 | -0.051 | -0.051 | -0.089 | -0.051 | -0.0513 |
| | $g_{11}(x)$ | -1.794E-9 | -1.434E-13 | -0.002 | -1.776E-8 | 7.777E-10 | 9.015E-14 | -5.866E-5 | -0.013 | 1.205E-4 | 5.185E-8 |
| Best | | 2994.4711 | 2994.4711 | 2994.750 | 2994.4716 | 2994.471 | 2994.4711 | 2994.4825 | 3000.981 | 2994.4671 | 2994.4711 |
| Mean | | 2994.4711 | 2994.4711 | NA | 2994.4717 | NA | 2994.4711 | 2999.6524 | 3007.1997 | 2994.4671 | 2994.4744 |
| worst | | 2994.4711 | 2994.4717 | NA | 2994.4732 | NA | 2994.4711 | 2999.6524 | 3009.9 | 4973.8644 | 2994.5056 |
| std | | 3.334E-7 | 2.976E-7 | NA | 4.251E-4 | NA | 1.239E-7 | 1.56 | 4.963 | 721.518 | 7.4E-3 |
| NFE | | 32000 | 26000 | NA | 30000 | 15000 | 22000 | 6300 | 250000 | 15000 | 15150 |

## 4.4. Sensitivity analysis

In this subsection, we will explore the influence of different values assigned to various parameters on the performance of the REAL. To demonstrate its effectiveness, we have selected one optimization problem from the CEC 2022 Test Suite [7]. This problem is derived by shifting and rotating the well-known Levy function. The Levy function is presented as follows:

$$f(x) = \sin^2(\pi w_1) + \sum_{i=1}^{D-1}(w_i - 1)^2 [1 + 10\sin^2(\pi w_i - 1)] + (w_D - 1)^2[1 + \sin^2(2\pi w_D)], \quad (11)$$

where $w_i = 1 + \frac{x_i - 1}{4}, \forall i = 1, \ldots, D$ and $D$ is the dimension of this function. In this subsection, we will take $D = 10$. The shifted and rotated Levy function to be used subsequently is as follows:

$$F(x) = f\left(M\left(\frac{5.12(x - 0.5)}{100}\right)\right) + F^*, \quad (12)$$

where $F^* = 900$ and $x \in [-100,100]^D$. The main properties of the shifted and rotated Levy function include multi-modal, non-separable, and huge number of local optima. The 3-D map for 2-D shifted and rotated Levy function is depicted in Figure 3. It is a challenge to search the global optima of this function with 10 dimensions.

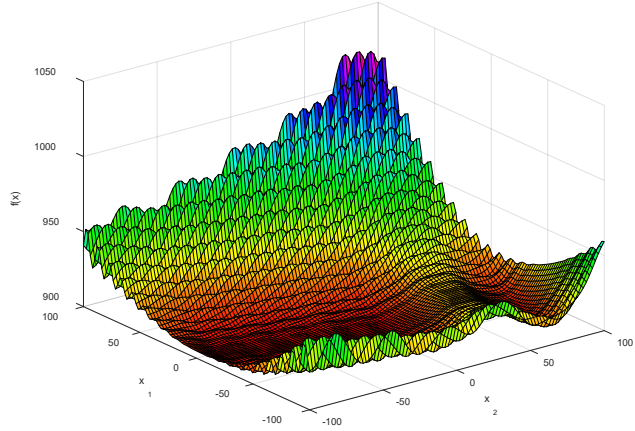

Fig. 3. 3-D map for 2-D shifted and rotated Levy function.

Initially, we address this problem using REAL with the given parameters: $n_X$=30, $\gamma$=6, T=1000, $n_{Rot}$=20, $n_{VS}$=10, $r_{A0} = 0.5$, $L_A^T = 1E - 40$, $r_P$=0.1, and $r_{Ex} = 0.5$. We set the termination condition such that when the absolute difference between $x_{vS}^{Best}$ and $x_{vS}^{Worst}$ is less than 1E-8 or the maximum iteration number is reached, REAL will terminate. In this case, REAL concludes its execution after 688 iterations, short of the maximum iteration limit of 1000. The achieved best objective function value is 900.0000000144273. The values of the corresponding ten elements of

the best solution $x$ are as follows: -61.232399714080586, -12.183568406870481, 50.259216843930304, 48.73364505462358, 1.788645752816961, -3.4670536406597536, 28.51212913481463, -3.742343827042961, -24.810503607058468, and 55.23631385323837. The required computational time is 1283 milliseconds.

The curves illustrating the changing values of $x_{vs}^{Best}$ and $x_{vs}^{Worst}$ over the 688 iterations are presented in Figure 4. It is evident that the associated values for these two indices decrease rapidly at the initial stages of the search process, which aligns with the typical behavior observed in other heuristic algorithms. Towards the end of the search, the values decrease more gradually, primarily due to the presence of large local optima near the global optimum.

To provide a comparative analysis of performance, Figure 5 displays the results for a standard agent, which holds a sequence of solutions generated successively by applying REAL's operations starting from an initial solution. The pronounced fluctuations in the objective function values associated with a typical agent highlight REAL's characteristic use of rotation, excursion, and perturbation operations to identify a promising solution. It's noteworthy that in Figure 5, the number of solutions linked to the agent is 890, exceeding the 688 iterations. This discrepancy arises because we only record the solution for the agent after each rotation or excursion operation, neglecting solutions resulting from perturbation.

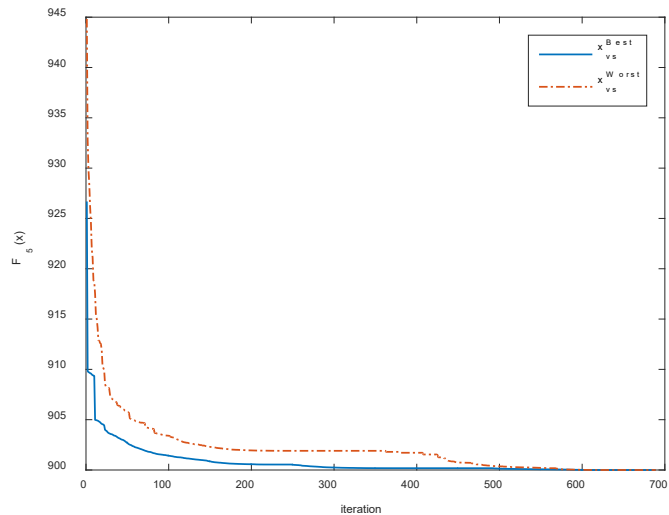

Fig. 4. The variation of objective function values associated with $x_{vs}^{Best}$ and $x_{vs}^{Worst}$.

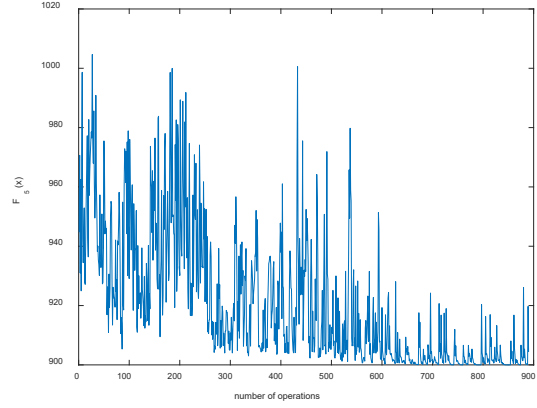

Fig. 5. The variation of objective values associated with an agent.

The points corresponding to the 890 solutions are illustrated in Figure 6, with a focus on the first three dimensions of the solution series. In Figure 6a, a 3-D map is presented, while the additional sub-figures (6b-6c) depict the projected scatter maps onto the $x_1 - x_2$, $x_1 - x_3$, and $x_2 - x_3$ planes, respectively. The maps in Figure 6 reveal that these solutions are widely dispersed across the entire search space, with many of them forming clusters around the final optima. This outcome aligns with the intended effect of REAL's implementation—thorough and random exploration of the search space, with a convergence towards an optimum at the end of the search.

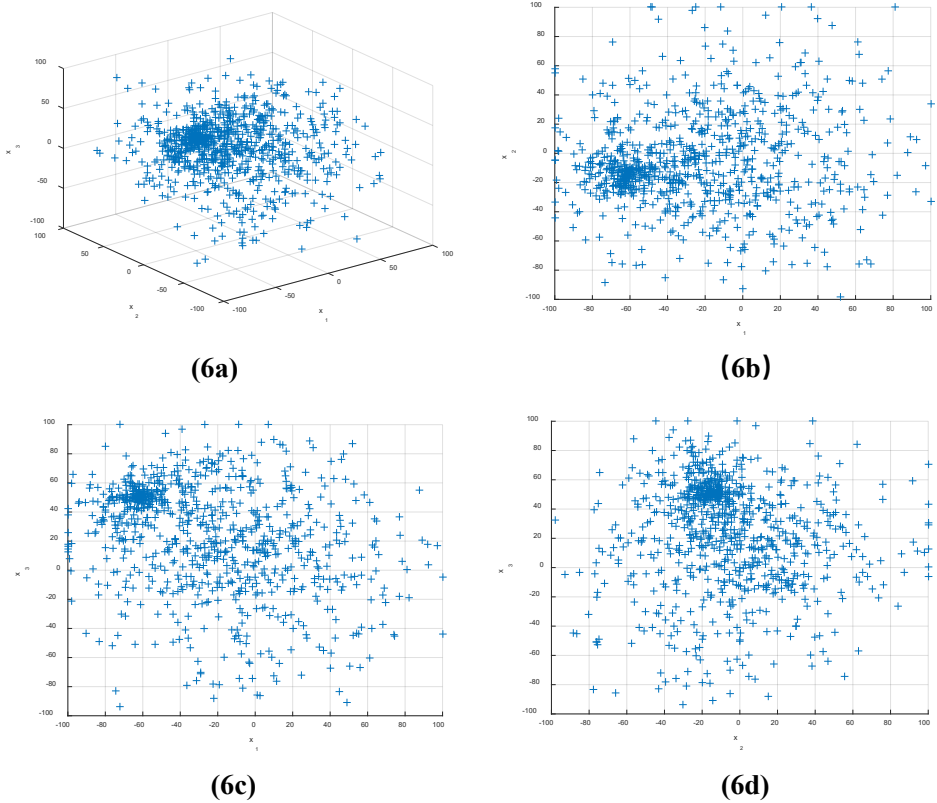

**(6a)** **(6b)**

**(6c)** **(6d)**

Fig. 6. The distribution of the 890 solutions of an agent $x$: (6a) $x_1$, $x_2$, and $x_3$ in 3-D; (6b) $x_1$-$x_2$ map; (6c) $x_1$-$x_3$ map; (6d) $x_2$-$x_3$ map.

Now, let's delve into the impact of parameters on REAL's performance. Nine parameters will be examined: $L_A^0$ (or $r_{A0}$), $L_A^T$ (or $r_{AT}$), $\gamma$, $r_P$, $r_{Ex}$, T, $n_{Rot}$, $n_{VS}$, and $n_X$. The relationships $L_A^0 = r_{A0}\|\beta - \alpha\|$ and $L_A^T = r_{AT}\|\beta - \alpha\|$ will guide our analysis. To simplify the study of $L_A^0$, we'll employ $r_{A0}$ due to the mentioned relationship. Similarly, for $L_A^T$, we'll utilize $r_{AT}$.

We will maintain the benchmark values for parameters in this section, with the exception of changing $T$ from 1000 to 500 to ensure diverse final results across different cases. With the specified parameter values, the Number of Function Evaluations (NFE) is approximately $2.5\, n_X\, T$=37,500 for a typical REAL run. Considering the feasible search space defined by $[-100, 100]^{10}$, it is straightforward to calculate $\|\beta - \alpha\| = \sqrt[2]{10 \times [100 - (-100)]^2} = 632.4555320336759$ for the 10-dimensional problem (12). The initial perturbation amplitude $L_A^0$ satisfies $r_{A0}\|\beta - \alpha\| = 0.5 \times \|\beta - \alpha\| \approx 316.227766$ based on the benchmark value of the parameter $r_{A0}$.

To streamline the table, the objective function value $F(x^*)$ of function (12) will be replaced by $F(x^*) - F^*$ in Tables 13-21. In simpler terms, the closer the values of the mean, best, and worst objective function values presented in Tables 13-21 are to zero, the better the performance.

Table 13 illustrates the REAL's performance with an increasing number of agents. As the number $n_x$ grows, the mean objective function value of 30 independent program runs for solving problem (12) gradually decreases. However, when it reaches the minimum at $n_x$=80, it starts to rise. This phenomenon occurs because the increasing number of agents can enhance the search capabilities of the REAL up to a certain point. Still, when it becomes too large, the REAL might prematurely converge to local optima, hindering further improvements in its search capacities.

**Table 13**

The performance of REAL with changing number of agents $n_X$.

| | \multicolumn{10}{c}{$n_X$} | | | | | | | | | |
|---|---|---|---|---|---|---|---|---|---|---|
| | 10 | 20 | 30 | 40 | 50 | 60 | 70 | 80 | 90 | 100 |
| Mean | 0.35046 | 0.25648 | 0.16048 | 0.09350 | 0.08236 | 0.06334 | 0.06063 | 0.05605 | 0.09195 | 0.07838 |
| Std. | 0.27676 | 0.36041 | 0.22858 | 0.16640 | 0.15221 | 0.14836 | 0.14496 | 0.11574 | 0.14752 | 0.14391 |
| Best | 0.00381 | 1.13E-12 | 7.95E-13 | 3.41E-13 | 2.27E-13 | 0.0 | 0.0 | 0.0 | 0.0 | 0.0 |
| Worst | 1.18176 | 1.45250 | 0.99817 | 0.54385 | 0.45446 | 0.54385 | 0.54385 | 0.45432 | 0.45432 | 0.54385 |
| NTM | 12500 | 25000 | 37500 | 50000 | 62500 | 75000 | 87500 | 100000 | 112500 | 125000 |

The value of $\gamma$ in equation (9) determines the shape of the learning efficiency during the growth period. A small value implies a wide span for the middle phase of rapidly accelerating learning efficiency. Conversely, a large $\gamma$ results in a narrower middle phase, indicating a more abrupt increase in learning efficiency over a shorter duration. The data in Table 14 reveal that when $\gamma$ equals 10, the mean objective function value is the smallest. Along with this lowest mean is the lowest best objective function value, representing an optimum. However, it is noteworthy that as the value of $\gamma$ increases to a relatively large value, such as 25, the associated mean does not strictly rise. This suggests that a relatively large $\gamma$ (or a narrow middle acceleration stage) may introduce some variability in the final stage of the search. In other words, better solutions might be discovered towards the end of the search due to an extended plateau phase where excursions are more likely to occur.

**Table 14**

The performance of REAL with changing radius of the truncated range given the maximal NFE=37,500.

| | \multicolumn{8}{c}{The radius of the truncated range $\gamma$} | | | | | | | |
|---|---|---|---|---|---|---|---|---|
| | 5 | 10 | 15 | 20 | 25 | 30 | 35 | 40 |
| Mean | 0.21471 | 0.07357 | 0.12561 | 0.18833 | 0.09959 | 0.17599 | 0.24405 | 0.15812 |
| Std. | 0.22192 | 0.07775 | 0.13849 | 0.19385 | 0.05576 | 0.18207 | 0.21512 | 0.17351 |
| Best | 7.95E-13 | 0.0 | 8.37E-4 | 0.00496 | 0.00848 | 0.01425 | 0.01180 | 0.00766 |
| Worst | 0.54515 | 0.27054 | 0.63678 | 0.72122 | 0.21730 | 0.74661 | 0.65030 | 0.73923 |

Table 15 demonstrates the impact of the maximum number of iterations on the REAL's performance. The results are unsurprising: as $T$ increases, the performance of the REAL gradually improves at the expense of increased computational time.

**Table 15**

The performance of REAL with changing maximal number of iterations $T$.

| | The maximal number of iterations $T$ | | | | | | | | | |
|---|---|---|---|---|---|---|---|---|---|---|
| | 100 | 200 | 300 | 400 | 500 | 600 | 700 | 800 | 900 | 1000 |
| Mean | 0.68710 | 0.32654 | 0.18257 | 0.18563 | 0.16048 | 0.15994 | 0.08666 | 0.10511 | 0.07028 | 0.03902 |
| Std. | 0.52069 | 0.30067 | 0.22563 | 0.22305 | 0.22858 | 0.22887 | 0.15823 | 0.19873 | 0.12265 | 0.08961 |
| Best | 0.10811 | 0.00652 | 2.61E-12 | 1.59E-12 | 7.95E-13 | 2.27E-13 | 1.13E-13 | 1.13E-13 | 0.0 | 0.0 |
| Worst | 2.19716 | 1.16944 | 0.73826 | 0.69559 | 0.99817 | 0.82505 | 0.54385 | 0.89493 | 0.45432 | 0.45432 |
| NFE | 7500 | 15000 | 22500 | 30000 | 37500 | 45000 | 52500 | 60000 | 67500 | 75000 |

The size of visible spots affects the REAL in two ways. On one hand, a relatively large size introduces diversity to the excursion operation, potentially leading to better solutions in the end. This diversity helps REAL avoid premature convergence to local optima. On the other hand, a large size implies increased computational cost and more dispersed search energy. The data in Table 16 reveals that when $n_{VS}$ is set to 3, the REAL achieves its best performance in terms of the mean, best, and worst objective function values. However, as the size increases to 11, the mean objective function value remains relatively low compared to others. This observation aligns with the earlier analysis of the influence of $n_{VS}$.

**Table 16**

The performance of REAL with changing size of the list of visible spots $n_{VS}$ with NFE=75000.

| | The size of the list of visible spots $n_{VS}$ | | | | | | | | | |
|---|---|---|---|---|---|---|---|---|---|---|
| | 3 | 4 | 5 | 6 | 7 | 8 | 9 | 10 | 11 | 12 |
| Mean | 0.03902 | 0.28097 | 0.24742 | 0.18224 | 0.20760 | 0.13943 | 0.17965 | 0.16048 | 0.07204 | 0.15426 |
| Std. | 0.08961 | 0.34057 | 0.27504 | 0.23056 | 0.20963 | 0.18356 | 0.21958 | 0.22858 | 0.13689 | 0.21315 |
| Best | 0.0 | 1.13E-13 | 2.27E-13 | 3.41E-13 | 5.68E-13 | 4.55E-13 | 3.41E-13 | 7.95E-13 | 4.55E-13 | 5.68E-13 |
| Worst | 0.45432 | 1.27516 | 0.99832 | 0.99818 | 0.71587 | 0.54386 | 0.71879 | 0.99817 | 0.54392 | 0.54405 |

The impact of the size of $M_{Rot}$ on the REAL's performance is illustrated in Table 17. When $n_{Rot}$ is set to 30, the mean achieves its lowest value. However, when $n_{Rot}$ is set to 45, the best objective function value reaches its lowest point. Additionally, it is noteworthy that when the size of $M_{Rot}$ is 15, both the mean and the best objective function values are the second lowest. This observation suggests that an optimal size for the rotation matrices should not be excessively large or small. In general, selecting a size from the interval [10,30] appears suitable for this problem.

**Table 17**

The performance of REAL with the changing size of the set of the orthogonal rotation matrices.

| | The size of the set of the orthogonal rotation matrices $n_{Rot}$ | | | | | | | | | |
|---|---|---|---|---|---|---|---|---|---|---|
| | 5 | 10 | 15 | 20 | 25 | 30 | 35 | 40 | 45 | 50 |
| Mean | 0.18572 | 0.12943 | 0.09877 | 0.16048 | 0.11191 | 0.08632 | 0.10272 | 0.11336 | 0.16506 | 0.13856 |
| Std. | 0.26544 | 0.19843 | 0.17599 | 0.22858 | 0.18378 | 0.15751 | 0.16500 | 0.18199 | 0.24577 | 0.18320 |
| Best | 9.09E-13 | 4.55E-13 | 2.27E-13 | 7.95E-13 | 4.54E-13 | 2.27E-13 | 5.68E-13 | 3.41E-13 | 1.14E-13 | 1.13E-12 |
| Worst | 1.17222 | 0.71587 | 0.80721 | 0.99817 | 0.54389 | 0.54386 | 0.54385 | 0.54385 | 0.99818 | 0.54745 |

The impact of the initial perturbation amplitude, $L_A^0$, is illustrated in Table 18. Mean objective values associated with relatively high $L_A^0$ are comparatively low, while the best objective value occurs at a low value of $L_A^0$ (or $r_{A0}$ equal to 0.2). If we consider the second-best values of the means and the bests, a relatively high value of $L_A^0$ proves to be a favorable choice. Based on this observation, it is evident that multiple trials are necessary to pinpoint an appropriate value for $L_A^0$.

**Table 18**

The performance of REAL with changing length of the initial perturbation amplitude $L_A^0$.

| | The size of $r_{A0}$ | | | | | | | | | |
|---|---|---|---|---|---|---|---|---|---|---|
| | 0.1 | 0.2 | 0.3 | 0.4 | 0.5 | 0.6 | 0.7 | 0.8 | 0.9 | 1.0 |
| Mean | 0.36885 | 0.25801 | 0.21756 | 0.13768 | 0.16048 | 0.07787 | 0.12132 | 0.10309 | 0.02018 | 0.11271 |
| Std. | 0.36193 | 0.28364 | 0.29067 | 0.20674 | 0.22858 | 0.13716 | 0.17995 | 0.19523 | 0.05687 | 0.14069 |
| Best | 7.71E-7 | 2.27E-13 | 5.68E-13 | 4.55E-13 | 7.95E-13 | 5.68E-13 | 7.95E-13 | 9.09E-13 | 3.41E-13 | 1.72E-12 |
| Worst | 1.17155 | 0.81755 | 1.25972 | 0.71623 | 0.99817 | 0.45432 | 0.54385 | 0.71587 | 0.24726 | 0.54385 |

The final perturbation amplitude, $L_A^T$, working in tandem with the initial amplitude, helps determine the actual amplitude during the search. The impact of different values of the final perturbation amplitude is presented in Table 19. The data reveal an intriguing trend regarding the influence of $L_A^T$. When its value is small, confined to the range $[1E-40, 1]$, the best objective function value is relatively small. However, when its value increases, as reflected in the value of $r_{AT}$, constrained to $[0.01, 0.2]$, the best objective function value becomes relatively large. The mean objective function values present a different perspective on the influence of $L_A^T$. When $r_{AT}$ equals 0.05, the associated mean is the smallest. In practical terms, if achieving a final best solution is the goal, a low value of $L_A^T$ and several independent runs of the REAL should be the optimal choice.

Table 19

The performance of REAL with changing final perturbation amplitude $L_A^T$ (or related parameter $r_{AT}$).

|  | $L_A^T$ | | | | | $r_{AT}$ | | | | |
| --- | --- | --- | --- | --- | --- | --- | --- | --- | --- | --- |
|  | 1E-40 | 1E-10 | 1E-2 | 1E-1 | 1E+0 | 0.01 | 0.05 | 0.1 | 0.15 | 0.2 |
| Mean | 0.16048 | 0.13621 | 0.09971 | 0.12559 | 0.13343 | 0.16605 | 0.07359 | 0.11166 | 0.14473 | 0.13563 |
| Std. | 0.22858 | 0.20712 | 0.16124 | 0.20087 | 0.19110 | 0.25467 | 0.11617 | 0.16857 | 0.21604 | 0.19416 |
| Best | 7.95E-13 | 2.27E-13 | 4.55E-13 | 5.68E-13 | 7.96E-13 | 1.71E-12 | 1.03E-11 | 6.45E-11 | 1.18E-11 | 3.01E-11 |
| Worst | 0.99817 | 0.71587 | 0.54385 | 0.54534 | 0.64264 | 0.93789 | 0.45432 | 0.54385 | 0.80540 | 0.54385 |

The perturbation rate, $r_P$, is employed to decide whether to perform a perturbation operation for each element of a feasible solution $x$. Its function is analogous to the mutation rate used in GA. Table 20 illustrates the impact of different values of $r_P$ on the performance of the REAL. When the value of $r_P$ is chosen from the interval $[0.15, 0.25]$, the REAL exhibits better performance in terms of the mean and best objective function values. It's important to note that the interval specified for this problem may not be suitable for other problems.

Table 20

The performance of REAL with the changing perturbation rate $r_P$.

|  | The perturbation rate $r_P$ | | | | | | | | | |
| --- | --- | --- | --- | --- | --- | --- | --- | --- | --- | --- |
|  | 0.001 | 0.005 | 0.01 | 0.05 | 0.1 | 0.15 | 0.2 | 0.25 | 0.3 | 0.4 |
| Mean | 0.54639 | 0.45587 | 0.23863 | 0.11756 | 0.16048 | 0.11488 | 0.11206 | 0.04696 | 0.10012 | 0.06269 |
| Std. | 0.51941 | 0.33941 | 0.26794 | 0.17545 | 0.22858 | 0.19596 | 0.16893 | 0.09347 | 0.16315 | 0.14429 |
| Best | 7.18E-4 | 5.28E-6 | 4.04E-11 | 5.68E-13 | 7.95E-13 | 1.14E-13 | 4.55E-13 | 4.55E-13 | 1.02E-12 | 1.48E-12 |
| Worst | 1.80844 | 1.04467 | 0.96169 | 0.54433 | 0.99817 | 0.80542 | 0.45432 | 0.49697 | 0.45432 | 0.67464 |

The excursion approaching rate, $r_{Ex}$, is utilized to regulate the amplitude of the excursion operation. A larger value of $r_{Ex}$ implies that the solution resulting from an excursion will be closer to the chosen visible spot. In Table 21, the impact of $r_{Ex}$ on the performance of the REAL is depicted. The data indicate that when the value of $r_{Ex}$ is relatively small, the corresponding mean objective function values are also relatively small. This suggests that, during the implementation of the REAL, a relatively small $r_{Ex}$ would be a preferable choice.

Table 21

The performance of REAL with the excursion approaching rate $r_{Ex}$.

|  | The excursion approaching rate $r_{Ex}$ | | | | | | | | | |
| --- | --- | --- | --- | --- | --- | --- | --- | --- | --- | --- |
|  | 0.1 | 0.2 | 0.3 | 0.4 | 0.5 | 0.6 | 0.7 | 0.8 | 0.9 | 1.0 |
| Mean | 0.17032 | 0.03965 | 0.03659 | 0.16325 | 0.16048 | 0.30961 | 0.36081 | 0.45145 | 0.65484 | 1.04438 |
| Std. | 0.21152 | 0.05592 | 0.05048 | 0.19334 | 0.22858 | 0.28244 | 0.27089 | 0.33760 | 0.45287 | 0.55134 |
| Best | 0.00151 | 5.68E-12 | 3.41E-13 | 3.41E-13 | 7.95E-13 | 2.27E-12 | 9.09E-13 | 2.27E-13 | 5.68E-13 | 0.17906 |
| Worst | 0.80835 | 0.18048 | 0.17967 | 0.54385 | 0.99817 | 1.19390 | 0.90865 | 1.17020 | 1.71490 | 2.50371 |

Now, let's consolidate the preceding investigation on the influence of parameters on the performance of the REAL. We can categorize them into three main groups.

In the first group, the optimal selection for a suitable value of parameters typically resides at the extremes of the specified range. This group encompasses $T$, $n_{VS}$, $n_x$, $\gamma$, and $r_{Ex}$. The trend is particularly pronounced in $T$ and $n_{VS}$ and less distinct in the other three parameters.

The second group consists of two parameters, $n_{Rot}$ and $r_P$. The most effective values for each parameter fall within a specific interval. Values that are excessively large or small, outside of this interval, do not contribute to improved performance in the REAL implementation.

The third group comprises two perturbation amplitudes, $L_A^T$ and $L_A^0$. For these parameters, there is no clearly defined interval or specific end of the examined range that ensures superior REAL performance. The determination of an optimal choice necessitates multiple trials involving different values of $L_A^T$ and $L_A^0$.

It is important to note that the sensitivity analysis conducted in this section pertains solely to the alteration of a single parameter while holding others constant. Observations may differ when two or more parameters change simultaneously. Given the insights gleaned from this analysis, a more comprehensive study in the future is warranted to elucidate the effects of various combinations of these parameters.

## 5. Conclusions

In this paper, we introduced a novel general-purpose heuristic algorithm named the Rotation Excursion Algorithm with Learning (REAL), motivated by the effectiveness of the generational mechanism observed in the CEC optimization problem series. The REAL incorporates three key operations to emulate this generational mechanism. The rotation operation manipulates solutions by rotating them using a pre-generated orthogonal rotation matrix. The excursion operation strategically positions solutions near historical high-performing solutions, known as visible spots. Additionally, the perturbation operation conducts neighborhood exploration for each solution. To optimize learning efficiency, REAL utilizes a variant of the Sigmoid function, determining whether to perform an excursion operation at each iteration and adjusting the perturbation amplitude as the search progresses.

Applying REAL to 67 benchmark optimization problems yielded promising results, surpassing other algorithms in numerous instances, with the highest Friedman mean Rank. When addressing real-life engineering optimization problems, REAL demonstrated superior performance in terms of mean, best, and worst final feasible solutions. A sensitivity analysis of the nine main parameters categorized them into three groups, providing valuable guidance for practitioners in selecting appropriate values.

While REAL exhibits promising performance, there are associated trade-offs. The matrix multiplication required for rotation operations, despite pre-prepared rotation matrices, incurs a time cost. Future approaches, such as the Levy random walk, may mitigate this. Another challenge lies in parameter assignment, especially for those with ambiguous optimal values, necessitating extensive trials.

Despite these considerations, we believe REAL can be a valuable tool for tackling complex optimization problems. Its versatility allows implementation across various domains, and its potential could be further harnessed by incorporating more intricate and problem-specific neighborhood search strategies beyond the current simple perturbation operation.

**Disclosure statement**

No potential conflict of interest was reported by the authors.

**Data availability statement**

The data that support the findings of this study are available from the corresponding author upon request.

# Appendix A (29 Standard Benchmark Optimization Problems):

**Table 1**

Unimodal benchmark functions.

| Functions | Dim | Range | $f_{min}$ |
|---|---|---|---|
| $f_1(x) = \sum_{i=1}^{n} x_i^2$ | 30 | $[-100,100]^n$ | $f(0,0,\ldots,0) = 0$ |
| $f_2(x) = \sum_{i=1}^{n} |x_i| + \prod_{i=1}^{n} |x_i|$ | 30 | $[-10,10]^n$ | $f(0,0,\ldots,0) = 0$ |
| $f_3(x) = \sum_{i=1}^{n} (\sum_{j=1}^{i} x_j)^2$ | 30 | $[-100,100]^n$ | $f(0,0,\ldots,0) = 0$ |
| $f_4(x) = \max_i \{|x_i|, 1 \leq i \leq n\}$ | 30 | $[-100,100]^n$ | $f(0,0,\ldots,0) = 0$ |
| $f_5(x) = \sum_{i=1}^{n-1}[100(x_{i+1} - x_i^2)^2 + (x_i - 1)^2]$ | 30 | $[-30,30]^n$ | $f(1,1,\ldots,1) = 0$ |
| $f_6(x) = \sum_{i=1}^{n}(\lfloor x_i + 0.5 \rfloor)^2$ | 30 | $[-100,100]^n$ | 0 |
| $f_7(x) = \sum_{i=1}^{n} i x_i^4 + random[0,1)$ | 30 | $[-1.28,1.28]^n$ | 0 |

**Table 2**

Multimodal benchmark functions.

| Functions | Dim | Range | $f_{min}$ |
|---|---|---|---|
| $f_8(x) = \sum_{i=1}^{n}[-x_i \sin(\sqrt{|x_i|})]$ | 30 | $[-500,500]^n$ | $-418.9829n$ |
| $f_9(x) = \sum_{i=1}^{n}[x_i^2 - 10\cos(2\pi x_i) + 10]$ | 30 | $[-5.12,5.12]^n$ | 0 |
| $f_{10}(x) = -20\exp\left(-0.2\sqrt{\frac{1}{n}\sum_{i=1}^{n} x_i^2}\right) - \exp\left(\frac{1}{n}\sum_{i=1}^{n} \cos(2\pi x_i)\right) + 20 + e$ | 30 | $[-32,32]^n$ | 0 |
| $f_{11}(x) = \sum_{i=1}^{d} \frac{x_i^2}{4000} - \prod_{i=1}^{d} \cos\left(\frac{x_i}{\sqrt{i}}\right) + 1$ | 30 | $[-512,512]^n$ | 0 |
| $f_{12}(x) = \frac{\pi}{n}\{10\sin^2(\pi y_1) + \sum_{i=1}^{n-1}(y_i - 1)^2[1 + 10\sin^2(\pi y_{i+1})] + (y_n - 1)^2\} + \sum_{i=1}^{n} u(x_i, 10, 100, 4)$ $y_i = 1 + \frac{x_i + 1}{4},\ u(x_i, a, k, m) = \begin{cases} k(x_i - a)^m & x_i > a \\ 0 & -a \leq x_i \leq a \\ k(-x_i - a)^m & x_i < -a \end{cases}$ | 30 | $[-50,50]^n$ | 0 |
| $f_{13}(x) = 0.1\{\sin^2(3\pi x_1) + \sum_{i=1}^{n}(x_i - 1)^2[1 + \sin^2(3\pi x_i + 1)] + (x_n - 1)^2[1 + \sin^2(2\pi x_n)]\} + \sum_{i=1}^{n} u(x_i, 5, 100, 4)$ | 30 | $[-50,50]^n$ | 0 |

**Table 3**

Fixed dimension multimodal benchmark functions.

| Functions | Dim | Range | $f_{min}$ |
|---|---|---|---|
| $f_{14}(x) = \left(\frac{1}{500} + \sum_{j=1}^{25} \frac{1}{j + \sum_{i=1}^{2}(x_i - a_{ij})^6}\right)^{-1}$ | 2 | $[-65.536, 65.536]^n$ | 1 |
| $f_{15}(x) = \sum_{i=1}^{11}\left[a_i - \frac{x_1(b_i^2 + b_i x_2)}{b_i^2 + b_i x_3 + x_4}\right]^2$ | 4 | $[-5,5]^n$ | 0.0003075 |
| $f_{16}(x) = 4x_1^2 - 2.1x_1^4 + x_1^6/3 + x_1 x_2 - 4x_2^2 + 4x_2^4$ | 2 | $[-5,5]^n$ | -1.0316285 |
| $f_{17}(x) = \left(x_2 - \frac{5.1}{4\pi^2}x_1^2 + \frac{5}{\pi}x_1 - 6\right)^2 + 10\left(1 - \frac{1}{8\pi}\right)\cos x_1 + 10$ | 2 | $[-5,5]^n$ | 0.398 |
| $f_{18}(x) = [1 + (x_1 + x_2 + 1)^2(19 - 14x_1 + 3x_1^2 - 14x_1 + 6x_1 x_2 + 3x_2^2)] \times$ $[30 + (2x_1 - 3x_2)^2(18 - 32x_1 + 12x_1^2 + 48x_2 - 36x_1 x_2 + 27x_2^2)]$ | 2 | $[-2,2]^n$ | 3 |
| $f_{19}(x) = -\sum_{i=1}^{4} c_i \exp\left(-\sum_{j=1}^{3} a_{ij}(x_j - p_{ij})^2\right)$ | 3 | $[0,1]^n$ | -3.86 |
| $f_{20}(x) = -\sum_{i=1}^{4} c_i \exp\left(-\sum_{j=1}^{6} a_{ij}(x_j - p_{ij})^2\right)$ | 6 | $[0,1]^n$ | -3.32 |
| $f_{21}(x) = -\sum_{i=1}^{5}[(x - a_i)(x - a_i)^T + c_i]^{-1}$ | 4 | $[0,10]^n$ | -10.1532 |
| $f_{22}(x) = -\sum_{i=1}^{7}[(x - a_i)(x - a_i)^T + c_i]^{-1}$ | 4 | $[0,10]^n$ | -10.4028 |
| $f_{23}(x) = -\sum_{i=1}^{10}[(x - a_i)(x - a_i)^T + c_i]^{-1}$ | 4 | $[0,10]^n$ | -10.5363 |

**Table 4**

Composite benchmark functions.

| Functions | Dim | Range | $f_{min}$ |
|---|---|---|---|
| $F_{24}$ (CF1) | 10 | [-5,5] | 0 |
| $f_1, f_2, f_3, \ldots, f_{10}$ =Sphere Function<br>$[\sigma_1, \sigma_2, \sigma_3, \ldots, \sigma_{10}] = [1,1,1,\ldots,1]$<br>$[\lambda_1, \lambda_2, \lambda_3, \ldots, \lambda_{10}] = [5/100, 5/100, 5/100, \ldots, 5/100]$ | | | |
| $F_{25}$ (CF2)<br>$f_1, f_2, f_3, \ldots, f_{10}$ =Griewank's Function<br>$[\sigma_1, \sigma_2, \sigma_3, \ldots, \sigma_{10}] = [1,1,1,\ldots,1]$<br>$[\lambda_1, \lambda_2, \lambda_3, \ldots, \lambda_{10}] = [5/100, 5/100, 5/100, \ldots, 5/100]$ | 10 | [-5,5] | 0 |
| $F_{26}$ (CF3)<br>$f_1, f_2, f_3, \ldots, f_{10}$ =Griewank's Function<br>$[\sigma_1, \sigma_2, \sigma_3, \ldots, \sigma_{10}] = [1,1,1,\ldots,1]$<br>$[\lambda_1, \lambda_2, \lambda_3, \ldots, \lambda_{10}] = [1,1,1,\ldots,1]$ | 10 | [-5,5] | 0 |
| $F_{27}$ (CF4)<br>$f_1, f_2$= Ackley's Function, $f_3, f_4$= Rastrigin's Function,<br>$f_5, f_6$= Weierstrass Function, $f_7, f_8$= Griewank's Function,<br>$f_9, f_{10}$= Sphere Function<br>$[\sigma_1, \sigma_2, \sigma_3, \ldots, \sigma_{10}] = [1,1,1,\ldots,1]$<br>$[\lambda_1, \lambda_2, \lambda_3, \ldots, \lambda_{10}] =[5/32, 5/32, 1, 1, 5/0.5, 5/0.5, 5/100, 5/100, 5/100, 5/100]$ | 10 | [-5,5] | 0 |
| $F_{28}$ (CF5)<br>$f_1, f_2$= Rastrigin s Function, $f_3, f_4$= Weierstrass Function,<br>$f_5, f_6$= Griewank s Function, $f_7, f_8$= Ackley s Function,<br>$f_9, f_{10}$ = Sphere Function<br>$[\sigma_1, \sigma_2, \sigma_3, \ldots, \sigma_{10}] = [1,1,1,\ldots,1]$<br>$[\lambda_1, \lambda_2, \lambda_3, \ldots, \lambda_{10}] =[1/5, 1/5, 5/0.5, 5/0.5, 5/100, 5/100, 5/32, 5/32, 5/100, 5/100]$ | 10 | [-5,5] | 0 |
| $f_{29}$ (CF6)<br>$f_1, f_2$= Rastrigin's Function, $f_3, f_4$= Weierstrass' Function,<br>$f_5, f_6$= Griewank's Function, $f_7, f_8$= Ackley's Function,<br>$f_9, f_{10}$= Sphere Function<br>$[\sigma_1, \sigma_2, \sigma_3, \ldots, \sigma_{10}] =[0.1, 0.2, 0.3, 0.4, 0.5, 0.6, 0.7, 0.8, 0.9, 1]$<br>$[\lambda_1, \lambda_2, \lambda_3, \ldots, \lambda_{10}] =[0.1*1/5, 0.2*1/5, 0.3*5/0.5, 0.4*5/0.5, 0.5*5/100,$<br>$0.6*5/100, 0.7*5/32, 0.8*5/32, 0.9*5/100, 1*5/100]$ | 10 | [-5,5] | 0 |

# Appendix B: (7 Classic Engineering Optimization Problems)

**B01. Multiple disc clutch brake design**

Consider variable $x = (r_i, r_o, t, F, Z)$

Minimize $f(x) = \pi(r_o^2 - r_i^2)t(Z + 1)\rho$

Subject to:

$g_1(x) = r_o - r_i - \Delta r \geq 0$, $g_2(x) = l_{max} - (Z + 1)(t + \delta) \geq 0$, $g_3(x) = p_{max} - p_{rz} \geq 0$,

$g_4(x) = p_{max}v_{srmax} - p_{rz}v_{sr} \geq 0$, $g_5(x) = v_{srmax} - v_{sr} \geq 0$, $g_6(x) = T_{max} - T \geq 0$,

$g_7(x) = M_h - sM_s \geq 0$, $g_8(x) = T \geq 0$,

where $M_h = \frac{2}{3}\mu F Z \frac{r_o^3 - r_i^3}{r_o^2 - r_i^2}$, $p_{rz} = \frac{F}{\pi(r_o^2 - r_i^2)}$, $M_h = \frac{2\pi n(r_o^3 - r_i^3)}{90(r_o^2 - r_i^2)}$, $T = \frac{I_z \pi n}{30(M_h + M_f)}$,

$\Delta r$ = 20mm, $t_{max}$=3mm, $t_{min}$=1.5mm, $l_{max}$ = 30mm, $Z_{max}$=10, $v_{srmax}$=10m/s, $\mu$ = 0.5,

s = 1.5, $M_s$=40Nm, $M_f$=3Nm, n=250rpm, $p_{max}$=1MPa, $I_z$=55kg mm$^2$, $T_{max}$=15s, $F_{max}$=1000N,

$r_{min}$=55mm, $r_{o\,max}$=110mm, $\rho$ = 7800kg/$m^3$.

Variable range $r_i \in \{60, 61, 62, \ldots, 80\}$, $r_o \in \{90, 91, 92, \ldots, 110\}$, $t \in \{1.0, 1.5, 2.0, 2.5, 3\}$, $F \in \{600, 610, 620, \ldots, 1000\}$,

$Z \in \{2,3,4,5,6,7,8,9\}$.

**B02. Robot gripper**

Minimize $f(x) = \max_z F_k(x, z) - \min_z F_k(x, z)$

Subject to:

$g_1(x) = Y_{min} - y(x, Z_{max}) \geq 0$, $g_2(x) = y(x, Z_{max}) \geq 0$, $g_3(x) = y(x, 0) - Y_{max} \geq 0$,

$g_4(x) = Y_G - y(x, 0) \geq 0$, $g_5(x) = (a + b)^2 - l^2 - e^2 \geq 0$,

$g_6(x) = (l - Z_{max})^2 + (a - e)^2 - b^2 \geq 0$, $g_7(x) = l - Z_{max} \geq 0$,

$g_8(x) = d(l, Z_{max}, e) + b - a \geq 0,$

$g_9(x) = d(l, 0, e) + b - a \geq 0,$

$g_{10}(x) = b + a - d(l, 0, e) \geq 0,$

where $d(l, z, e) = \sqrt{(l-z)^2 + e^2}$, $\emptyset = \tan^{-1}\left(\frac{e}{l-z}\right)$, $\alpha = \cos^{-1}\left(\frac{a^2+d^2-b^2}{2ad}\right) + \emptyset$,

$\beta = \cos^{-1}\left(\frac{b^2+d^2-a^2}{2bd}\right) - \emptyset$, $F_k = \frac{Pb\sin(\alpha+\beta)}{2c\cos\alpha}$, $y(x, z) = 2(e + f + c\sin(\beta + \delta))$,

$Y_{min} = 50$, $Y_{max} = 100$, $Y_G = 150$, $Z_{max} = 100$, $P = 100$,

$10 \leq a, b, f \leq 150$, $100 \leq c \leq 200$, $0 \leq e \leq 50$, $100 \leq l \leq 300$, $1 \leq \delta \leq 3.14$.

The constraints $g_8(x)$, $g_9(x)$ and $g_{10}(x)$ are newly added to guarantee the realizability of the results.

**B03. Rolling element bearing**

Maximize $C_d = f_c Z^{2/3} D_b^{1.8}$  if $D_b \leq 25.4$mm

$\quad\quad\quad\quad C_d = 3.64 f_c Z^{2/3} D_b^{1.4}$  if $D_b > 25.4$mm

Subject to:

$g_1(x) = \frac{\emptyset_o}{2\sin^{-1}(D_b/D_m)} - Z + 1 \geq 0$, $g_2(x) = 2D_b - K_{Dmin}(D - d) \geq 0$,

$g_3(x) = K_{Dmax}(D - d) - 2D_b \geq 0$, $g_4(x) = \zeta B_w - D_b \leq 0$,

$g_5(x) = D_m - 0.5(D + d) \geq 0$, $g_6(x) = (0.5 + e)(D + d) - D_m \geq 0$,

$g_7(x) = 0.5(D - D_m - D_b) - \varepsilon D_b \geq 0$, $g_8(x) = f_i \geq 0.515$, $g_9(x) = f_o \geq 0.515$,

where

$f_c = 37.91 \left\{ 1 + \left[ 1.04 \left(\frac{1-\gamma}{1+\gamma}\right)^{1.72} \left(\frac{f_i(2f_o-1)}{f_o(2f_i-1)}\right)^{0.41} \right]^{10/3} \right\}^{-0.3} \left(\frac{\gamma^{0.3}(1-\gamma)^{1.39}}{(1+\gamma)^{1/3}}\right) \left(\frac{2f_i}{2f_i-1}\right)^{0.41}$,

$\phi_o = 2\pi - 2\cos^{-1} \frac{\left[\frac{(D-d)}{2} - 3\left(\frac{T}{4}\right)\right]^2 + \left(\frac{D}{2} - \frac{T}{4} - D_b\right)^2 - \left(\frac{d}{2} + \frac{T}{4}\right)^2}{2\left[\frac{(D-d)}{2} - 3\left(\frac{T}{4}\right)\right]\left(\frac{D}{2} - \frac{T}{4} - D_b\right)}$,

$\gamma = \frac{D_b \cos\alpha}{D_m}$, $f_i = \frac{r_i}{D_b}$, $f_o = \frac{r_o}{D_b}$, $T = D - d - 2D_b$, $D = 160$, $d = 90$, $B_w = 30$, $\alpha = 0$,

$0.5(D + d) \leq D_m \leq 0.6(D + d)$, $0.15(D - d) \leq D_b \leq 0.45(D - d)$, $Z \in \{4, 5, 6, \ldots, 50\}$,

$0.515 \leq f_i, f_o \leq 0.6$, $0.4 \leq K_{Dmin} \leq 0.5$, $0.6 \leq K_{Dmax} \leq 0.7$, $0.3 \leq \varepsilon \leq 0.4$,

$0.02 \leq e \leq 0.1$, $0.6 \leq \zeta \leq 0.85$.

**B04. Hydrodynamic thrust bearing**

Minimize $f(x) = \frac{QP_0}{0.7} + E_f$

Subject to: $g_1(x) = W - W_s \geq 0$, $g_2(x) = P_{max} - P_0 \geq 0$, $g_3(x) = \Delta T_{max} - \Delta T \geq 0$, $g_4(x) = h - h_{min} \geq 0$, $g_5(x) = R - R_0 \geq 0$, $g_6(x) = 0.001 - \frac{\gamma}{gP_0}\left(\frac{Q}{2\pi Rh}\right)^2 \geq 0$, $g_7(x) = 5000 - \frac{W}{\pi(R^2 - R_0^2)} \geq 0$,

Where $\gamma = 0.0307$, $C = 0.5$, $n = -3.55$, $C_1 = 10.04$, $W_s = 101000$, $P_{max} = 1000$, $\Delta T_{max} = 50$, $h_{min} = 0.001$,

$g = 386.4$, $N = 750$, $1 \leq R, R_0, Q \leq 16$, $1e-6 \leq \mu \leq 16e-6$.

**B05. Belleville spring**

Minimize $f(x) = 0.07075\pi(D_e^2 - D_i^2)t$

Subject to:

$g_1(x) = S - \frac{4E\delta_{max}}{(1-\mu^2)\alpha D_e^2}\left[\beta\left(h - \frac{\delta_{max}}{2}\right) + \gamma t\right] \geq 0$,

$g_2(x) = \frac{4E\delta_{max}}{(1-\mu^2)\alpha D_e^2}\left[\left(h - \frac{\delta_{max}}{2}\right)(h - \delta_{max})t + t^3\right] - P_{max} \geq 0$,

$g_3(x) = \delta_l - \delta_{max} \geq 0$, $g_4(x) = H - h - t \geq 0$, $g_5(x) = D_{max} - D_e \geq 0$,

$g_6(x) = D_e - D_i \geq 0$, $g_7(x) = 0.3 - \frac{1}{D_e - D_i} \geq 0$,

where $\alpha = \frac{6}{\pi \ln K}\left(\frac{K-1}{K}\right)^2$, $\beta = \frac{6}{\pi \ln K}\left(\frac{K-1}{\ln K} - 1\right)$, $\gamma = \frac{6}{\pi \ln K}\left(\frac{K-1}{2}\right)$,

$P_{max} = 5400$lb, $\delta_{max} = 0.2$in., $S=200$kPsi, $E=30e6$psi, $\mu = 0.3$, $H=0.2$in., $D_{max}=12.01$in.,

$K = \frac{D_e}{D_i}$, $\delta_l = f(a)h$, $a = h/t$. $0.01 \leq t \leq 6.0$, $0.05 \leq h \leq 0.5$, $5.0 \leq D_i, D_e \leq 15.0$.

Values of $f(a)$ vary as shown in Table 7.

Table 7

Variation of $f(a)$ with $a$.

| $a$ | $\leq 1.4$ | 1.5 | 1.6 | 1.7 | 1.8 | 1.9 | 2.0 | 2.1 | 2.2 | 2.3 | 2.4 | 2.5 | 2.6 | 2.7 | $\geq 2.8$ |
|---|---|---|---|---|---|---|---|---|---|---|---|---|---|---|---|
| $f(a)$ | 1 | 0.85 | 0.77 | 0.71 | 0.66 | 0.63 | 0.6 | 0.58 | 0.56 | 0.55 | 0.53 | 0.52 | 0.51 | 0.51 | 0.5 |

**B06. Step-cone pulley**

Minimize $f(x) = \rho w \left\{ d_1^2 \left[1 + \left(\frac{N_1}{N}\right)^2\right] + d_2^2 \left[1 + \left(\frac{N_2}{N}\right)^2\right] + d_3^2 \left[1 + \left(\frac{N_3}{N}\right)^2\right] + d_4^2 \left[1 + \left(\frac{N_4}{N}\right)^2\right] \right\}$

Subject to:

$h_1(x) = C_1 - C_2 = 0$, $h_2(x) = C_1 - C_3 = 0$, $h_3(x) = C_1 - C_4 = 0$,

$g_{1,2,3,4}(x) = R_i \geq 2$, $g_{5,6,7,8}(x) = P_i \geq (0.75 * 745.6998)$,

where $C_i$ indicates the length of the belt to obtain speed $N_i$ and is given by

$C_i = \frac{\pi d_i}{2}\left(1 + \frac{N_i}{N}\right) + \frac{\left(\frac{N_i}{N}-1\right)^2 d_i^2}{4a} + 2a$, $i = (1,2,3,4)$.

$R_i$ is the tension ratio and is given by

$R_i = \exp\left(\mu \left[\pi - 2\sin^{-1}\left(\left(\frac{N_i}{N} - 1\right)\frac{d_i}{2a}\right)\right]\right)$, $i = (1,2,3,4)$.

$P_i = stw\left[1 - exp\left[-\mu\left\{\pi - 2\sin^{-1}\left(\left(\frac{N_i}{N} - 1\right)\frac{d_i}{2a}\right)\right\}\right]\right]\frac{\pi d_i N_i}{60}$, $i = (1,2,3,4)$.

$\rho = 7200$kg/$m^3$, $a = 3$m, $\mu = 0.35$, $s = 1.75$MPa, $t = 8$mm.

$1$mm $\leq d_1, d_2, d_3, d_4, w \leq 100$mm.

**B07. Speed reducer design**

Minimize $f(x) = 0.7854 x_1 x_2^2 (3.3333 x_3^2 + 14.9334 x_3 - 43.0934) - 1.508 x_1 (x_6^2 + x_7^2)$

$+ 7.4777(x_6^3 + x_7^3) + 0.7854(x_4 x_6^2 + x_5 x_7^2)$

Subject to:

$g_1(x) = \frac{27}{x_1 x_2^2 x_3} - 1 \leq 0$, $g_2(x) = \frac{397.5}{x_1 x_2^2 x_3^2} - 1 \leq 0$, $g_3(x) = \frac{1.93 x_4^3}{x_2 x_3 x_6^4} - 1 \leq 0$, $g_4(x) = \frac{1.93 x_5^3}{x_2 x_3 x_7^4} - 1 \leq 0$,

$g_5(x) = \frac{\sqrt{[745 x_4/(x_2 x_3)]^2 + 16.9e6}}{110 x_6^3} - 1 \leq 0$, $g_6(x) = \frac{\sqrt{[745 x_5/(x_2 x_3)]^2 + 157.5e6}}{85 x_7^3} - 1 \leq 0$,

$g_7(x) = \frac{x_2 x_3}{40} - 1 \leq 0$, $g_8(x) = \frac{5 x_2}{x_1} - 1 \leq 0$, $g_9(x) = \frac{x_1}{12 x_2} - 1 \leq 0$,

$g_{10}(x) = \frac{1.5 x_6 + 1.9}{x_4} - 1 \leq 0$, $g_{11}(x) = \frac{1.1 x_7 + 1.9}{x_5} - 1 \leq 0$,

where $2.6 \leq x_1 \leq 3.6$; $0.7 \leq x_2 \leq 0.8$; $17 \leq x_3 \leq 28$; $7.3 \leq x_4 \leq 8.3$; $7.8 \leq x_5 \leq 8.3$;

$2.9 \leq x_6 \leq 3.9$; $5.0 \leq x_7 \leq 5.5$.

# References:

[1] S.H. Zanakis, J.R. Evans, Heuristic "Optimization": Why, When, and How to Use It, Interfaces, 11

(1981) 84-91.

[2] D.H. Wolpert, W.G. Macready, No free lunch theorems for optimization, IEEE transactions on evolutionary computation, 1 (1997) 67-82.

[3] A.P. Piotrowski, J.J. Napiorkowski, A.E. Piotrowska, Choice of benchmark optimization problems does matter, Swarm and Evolutionary Computation, 83 (2023) 101378.

[4] J. Liang, B. Qu, P. Suganthan, Problem definitions and evaluation criteria for the CEC 2014 special session and competition on single objective real-parameter numerical optimization, 2013.

[5] N. Awad, M. Ali, J. Liang, B. Qu, P. Suganthan, Problem definitions and evaluation criteria for the CEC 2017 special session and competition on single objective bound constrained real-parameter numerical optimization, in: Technical Report, Nanyang Technological University Singapore, 2016, pp. 1-34.

[6] C. Yue, K.V. Price, P.N. Suganthan, J. Liang, M.Z. Ali, B. Qu, N.H. Awad, P.P. Biswas, Problem definitions and evaluation criteria for the CEC 2020 special session and competition on single objective bound constrained numerical optimization, Comput. Intell. Lab., Zhengzhou Univ., Zhengzhou, China, Tech. Rep, 201911 (2019).

[7] A. Kumar, K.V. Price, A.W. Mohamed, A.A. Hadi, P.N. Suganthan, Problem Definitions and Evaluation Criteria for the 2022 Special Session and Competition on Single Objective Bound Constrained Numerical Optimization, in, Nanyang Technological University, Singapore, 2021, pp. 1-20.

[8] J.H. Holland, Genetic algorithms, Scientific american, 267 (1992) 66-73.

[9] J. Kennedy, R. Eberhart, Particle swarm optimization, in: Proceedings of ICNN'95 - International Conference on Neural Networks, 1995, pp. 1942-1948 vol.1944.

[10] R. Storn, K. Price, Differential evolution–a simple and efficient heuristic for global optimization over continuous spaces, Journal of global optimization, 11 (1997) 341-359.

[11] S. He, E. Prempain, Q.H. Wu, An improved particle swarm optimizer for mechanical design optimization problems, Engineering Optimization, 36 (2004) 585-605.

[12] E. Mezura-Montes, C.A.C. Coello, A simple multimembered evolution strategy to solve constrained optimization problems, IEEE Transactions on Evolutionary Computation, 9 (2005) 1-17.

[13] D. Karaboğa, AN IDEA BASED ON HONEY BEE SWARM FOR NUMERICAL OPTIMIZATION, in, 2005.

[14] M. Dorigo, M. Birattari, T. Stutzle, Ant colony optimization, IEEE computational intelligence magazine, 1 (2006) 28-39.

[15] D. Simon, Biogeography-based optimization, IEEE transactions on evolutionary computation, 12 (2008) 702-713.

[16] X.-S. Yang, S. Deb, Cuckoo search via Lévy flights, in: 2009 World congress on nature & biologically inspired computing (NaBIC), Ieee, 2009, pp. 210-214.

[17] R. Oftadeh, M.J. Mahjoob, M. Shariatpanahi, A novel meta-heuristic optimization algorithm inspired by group hunting of animals: Hunting search, Computers & Mathematics with Applications, 60 (2010) 2087-2098.

[18] X.-S. Yang, Firefly algorithm, stochastic test functions and design optimisation, in: arXiv preprint arXiv:1003.1409, 2010, pp. 1-12.

[19] X.-S. Yang, A new metaheuristic bat-inspired algorithm, in: Nature inspired cooperative strategies for optimization (NICSO 2010), Springer, 2010, pp. 65-74.

[20] X.S. Yang, Nature-inspired Metaheuristic Algorithms, Luniver Press, 2010.

[21] X.-S. Yang, Flower Pollination Algorithm for Global Optimization, in: International Conference

on Unconventional Computation and Natural Computation, 2012.

[22] X.S. Yang, A. Hossein Gandomi, Bat algorithm: a novel approach for global engineering optimization, Engineering Computations, 29 (2012) 464-483.

[23] S. Mirjalili, S.M. Mirjalili, A. Lewis, Grey Wolf Optimizer, Advances in Engineering Software, 69 (2014) 46-61.

[24] S. Mirjalili, Moth-flame optimization algorithm: A novel nature-inspired heuristic paradigm, Knowledge-Based Systems, 89 (2015) 228-249.

[25] S. Mirjalili, A. Lewis, The Whale Optimization Algorithm, Advances in Engineering Software, 95 (2016) 51-67.

[26] S. Mirjalili, A.H. Gandomi, S.Z. Mirjalili, S. Saremi, H. Faris, S.M. Mirjalili, Salp Swarm Algorithm: A bio-inspired optimizer for engineering design problems, Advances in Engineering Software, 114 (2017) 163-191.

[27] G. Azizyan, F. Miarnaeimi, M. Rashki, N. Shabakhty, Flying Squirrel Optimizer (FSO): A novel SI-based optimization algorithm for engineering problems, Iranian Journal of Optimization, 11 (2019) 177-205.

[28] A.A. Heidari, S. Mirjalili, H. Faris, I. Aljarah, M. Mafarja, H. Chen, Harris hawks optimization: Algorithm and applications, Future Generation Computer Systems, 97 (2019) 849-872.

[29] B. Abdollahzadeh, F.S. Gharehchopogh, S. Mirjalili, African vultures optimization algorithm: A new nature-inspired metaheuristic algorithm for global optimization problems, Computers & Industrial Engineering, 158 (2021) 107408.

[30] W. Zhao, L. Wang, S. Mirjalili, Artificial hummingbird algorithm: A new bio-inspired optimizer with its engineering applications, Computer Methods in Applied Mechanics and Engineering, 388 (2022) 114194.

[31] F.A. Hashim, A.G. Hussien, Snake Optimizer: A novel meta-heuristic optimization algorithm, Knowledge-Based Systems, 242 (2022) 108320.

[32] T. Ray, K.M. Liew, Society and civilization: An optimization algorithm based on the simulation of social behavior, IEEE Transactions on Evolutionary Computation, 7 (2003) 386-396.

[33] K.S. Lee, Z.W. Geem, A new meta-heuristic algorithm for continuous engineering optimization: harmony search theory and practice, Computer Methods in Applied Mechanics and Engineering, 194 (2005) 3902-3933.

[34] R.V. Rao, V.J. Savsani, D.P. Vakharia, Teaching–learning-based optimization: A novel method for constrained mechanical design optimization problems, Computer-Aided Design, 43 (2011) 303-315.

[35] S.-X. He, Truss optimization with frequency constraints using the medalist learning algorithm, Structures, 55 (2023) 1-15.

[36] S.-X. He, Y.-T. Cui, Medalist learning algorithm for configuration optimization of trusses, Applied Soft Computing, 148 (2023) 110889.

[37] S.-X. He, Y.-T. Cui, A novel variational inequality approach for modeling the optimal equilibrium in multi-tiered supply chain networks, Supply Chain Analytics, 4 (2023) 100039.

[38] S.-X. He, Y.-T. Cui, Multiscale medalist learning algorithm and its application in engineering, Acta Mechanica, (2023).

[39] S. Kirkpatrick, C.D. Gelatt, M.P. Vecchi, Simulated annealing, Science, 220 (1983) 671-680.

[40] R. Formato, Central force optimization: a new metaheuristic with applications in applied electromagnetics, Progress In Electromagnetics Research, 77 (2007) 425-491.

[41] E. Rashedi, H. Nezamabadi-pour, S. Saryazdi, GSA: A Gravitational Search Algorithm, Information

Sciences, 179 (2009) 2232-2248.

[42] A. Kaveh, S. Talatahari, A novel heuristic optimization method: charged system search, Acta Mechanica, 213 (2010) 267-289.

[43] H. Eskandar, A. Sadollah, A. Bahreininejad, M. Hamdi, Water cycle algorithm – A novel metaheuristic optimization method for solving constrained engineering optimization problems, Computers & Structures, 110-111 (2012) 151-166.

[44] F.F. Moghaddam, R.F. Moghaddam, M. Cheriet, Curved space optimization: a random search based on general relativity theory, arXiv preprint arXiv, 1 (2012) 12082214.

[45] A. Hatamlou, Black hole: A new heuristic optimization approach for data clustering, Information sciences, 222 (2013) 175-184.

[46] A. Sadollah, A. Bahreininejad, H. Eskandar, M. Hamdi, Mine blast algorithm: A new population based algorithm for solving constrained engineering optimization problems, Applied Soft Computing, 13 (2013) 2592-2612.

[47] H. Salimi, Stochastic Fractal Search: A powerful metaheuristic algorithm, Knowledge-Based Systems, 75 (2015) 1-18.

[48] P. Savsani, V. Savsani, Passing vehicle search (PVS): A novel metaheuristic algorithm, Applied Mathematical Modelling, 40 (2016) 3951-3978.

[49] A. Kaveh, A. Dadras, A novel meta-heuristic optimization algorithm: Thermal exchange optimization, Advances in Engineering Software, 110 (2017) 69-84.

[50] A. Faramarzi, M. Heidarinejad, B. Stephens, S. Mirjalili, Equilibrium optimizer: A novel optimization algorithm, Knowledge-Based Systems, 191 (2020) 105190.

[51] W. Zhao, L. Wang, Z. Zhang, Artificial ecosystem-based optimization: a novel nature-inspired meta-heuristic algorithm, Neural Computing and Applications, 32 (2020) 9383-9425.

[52] H. Ma, H. Wei, Y. Tian, R. Cheng, X. Zhang, A multi-stage evolutionary algorithm for multi-objective optimization with complex constraints, Information Sciences, 560 (2021) 68-91.

[53] B.S. Yildiz, N. Pholdee, S. Bureerat, A.R. Yildiz, S.M. Sait, Enhanced grasshopper optimization algorithm using elite opposition-based learning for solving real-world engineering problems, Engineering with Computers, 38 (2022) 4207-4219.

[54] Y. Zhang, Elite archives-driven particle swarm optimization for large scale numerical optimization and its engineering applications, Swarm and Evolutionary Computation, 76 (2023) 101212.

[55] M. Ghasemi, M. Zare, A. Zahedi, M.-A. Akbari, S. Mirjalili, L. Abualigah, Geyser Inspired Algorithm: A New Geological-inspired Meta-heuristic for Real-parameter and Constrained Engineering Optimization, Journal of Bionic Engineering, (2023).

[56] F. Rezaei, H.R. Safavi, M. Abd Elaziz, S. Mirjalili, GMO: geometric mean optimizer for solving engineering problems, Soft Computing, 27 (2023) 10571-10606.

[57] K.O. Mohammed Aarif, P. Sivakumar, M.Y. Caffiyar, B.A.M. Hashim, C.M. Hashim, C.A. Rahman, Nature-Inspired Optimization Algorithms: Past to Present, in: J. Nayak, A.K. Das, B. Naik, S.K. Meher, S. Brahnam (Eds.) Nature-Inspired Optimization Methodologies in Biomedical and Healthcare, Springer International Publishing, Cham, 2023, pp. 1-32.

[58] X.-S. Yang, Nature-Inspired Algorithms in Optimization: Introduction, Hybridization, and Insights, in: X.-S. Yang (Ed.) Benchmarks and Hybrid Algorithms in Optimization and Applications, Springer Nature Singapore, Singapore, 2023, pp. 1-17.

[59] W. Korani, M. Mouhoub, Review on Nature-Inspired Algorithms, Operations Research Forum, 2 (2021) 36.

[60] R. Salomon, Re-evaluating genetic algorithm performance under coordinate rotation of benchmark functions. A survey of some theoretical and practical aspects of genetic algorithms, Biosystems, 39 (1996) 263-278.

[61] N. Hansen, S.D. Müller, P. Koumoutsakos, Reducing the time complexity of the derandomized evolution strategy with covariance matrix adaptation (CMA-ES), Evolutionary computation, 11 (2003) 1-18.

[62] R. Tanabe, A. Fukunaga, Success-history based parameter adaptation for Differential Evolution, in: 2013 IEEE Congress on Evolutionary Computation, 2013, pp. 71-78.

[63] A.W. Mohamed, A.A. Hadi, A.M. Fattouh, K.M. Jambi, LSHADE with semi-parameter adaptation hybrid with CMA-ES for solving CEC 2017 benchmark problems, in: 2017 IEEE Congress on Evolutionary Computation (CEC), 2017, pp. 145-152.

[64] S. Salkuti, B. Panigrahi, R. Kundu, R. Mukherjee, S. Debchoudhury, Energy and spinning reserve scheduling for a wind-thermal power system using CMA-ES with mean learning technique, International Journal of Electrical Power & Energy Systems, 53 (2013) 113–122.

[65] R. Tanabe, A. Fukunaga, Success-history based parameter adaptation for Differential Evolution, 2013.

[66] B.S. Yildiz, N. Pholdee, S. Bureerat, A.R. Yildiz, S.M. Sait, Robust design of a robot gripper mechanism using new hybrid grasshopper optimization algorithm, Expert Systems, 38 (2021) e12666.

[67] S. Krenich, Optimal Design of Robot Gripper Mechanism Using Force and Displacement Transmission Ratio, Applied Mechanics and Materials, 613 (2014) 117-125.

[68] C.A.C. Coello, TREATING CONSTRAINTS AS OBJECTIVES FOR SINGLE-OBJECTIVE EVOLUTIONARY OPTIMIZATION, Engineering Optimization, 32 (2000) 275-308.

[69] K. Deb, M. Goyal, Optimizing Engineering Designs Using a Combined Genetic Search, in: International Conference on Genetic Algorithms, 1997.

[70] J.N. Siddall, Optimal Engineering Design: Principles and Applications, Marcel Dekker, Inc., 1982.

[71] W. Gong, Z. Cai, D. Liang, Engineering optimization by means of an improved constrained differential evolution, Computer Methods in Applied Mechanics and Engineering, 268 (2014) 884-904.

[72] A.H. Gandomi, X.-S. Yang, A.H. Alavi, Cuckoo search algorithm: a metaheuristic approach to solve structural optimization problems, Engineering with Computers, 29 (2013) 17-35.

[73] B. Akay, D. Karaboga, Artificial bee colony algorithm for large-scale problems and engineering design optimization, Journal of Intelligent Manufacturing, 23 (2012) 1001-1014.